\documentclass[reqno,11 pt]{amsart}
\usepackage{graphicx} 
\usepackage{amsmath,amssymb,latexsym,float}
\usepackage{amsthm}
\usepackage{color}  
\usepackage[all]{xy}
\usepackage{tikz-cd}
\usepackage{multirow}
\usepackage{multicol}
\usepackage{longtable}
\usepackage{caption}
\usepackage{cite}
\usepackage{enumitem}
\usepackage{listings}
\usepackage{hyperref}
\usetikzlibrary{cd}
\usepackage[left=0.25in,right=0.25in, top= -0.5in, bottom= -0.5 in]{geometry}
\usepackage{array}
\usepackage{changepage}

\newcommand{\F}{{\mathbb{F}}}
\newcommand{\Z}{{\mathbb{Z}}}
\newcommand{\Q}{{\mathbb{Q}}} 
 
\newcommand{\QQ}{\overline{{\mathbb{Q}}}}

\newcommand{\p}{{\mathfrak{p}}}
\newcommand{\OO}{{\mathcal{O}}}

\newcommand{\Gal}{\mathrm{Gal}}
\newcommand{\Hom}{\mathrm{Hom}}
\newcommand{\Sel}{\mathrm{Sel}}

\newcommand{\rk}{\mathrm{rk} }

\newtheorem{lthm}{Theorem}
\newtheorem{lpro}{Proposition}
\newtheorem{lcor}{Corollary}

\makeatletter
\newcommand*{\rom}[1]{\expandafter\@slowromancap\romannumeral #1@}
\makeatother

\DeclareFontFamily{U}{wncy}{}
    \DeclareFontShape{U}{wncy}{m}{n}{<->wncyr10}{}
    \DeclareSymbolFont{mcy}{U}{wncy}{m}{n}
    \DeclareMathSymbol{\Sh}{\mathord}{mcy}{"58}
    		
\theoremstyle{plain}
\newtheorem{theorem}{Theorem}[section]
\newtheorem*{theorem*}{Theorem}
\newtheorem{proposition}[theorem]{Proposition}
\newtheorem{rem}[theorem]{Remark}
\newtheorem{lemma}[theorem]{Lemma}
\newtheorem{corollary}[theorem]{Corollary}

\newtheorem{conj}[theorem]{Conjecture}

\title{Relative $p$-class groups and $p$-Selmer groups }
\author[D.~De]{Debajyoti De}
\address[De]{Department of Mathematics, IIT Madras, India}
\email{debajyotide20@gmail.com}

\author[D.~Majumdar]{Dipramit Majumdar}
\address[Majumdar]{Department of Mathematics, IIT Madras, India}
\email{dipramit@iitm.ac.in}

\author[S.~Mondal]{Sudipa Mondal}
\address[Mondal]{Department of Mathematics, IIT Madras, India}
\email{sudipa.mondal123@gmail.com}


\begin{document}

\begin{abstract}
 Let $E$ be an elliptic curve with $j$-invariant $0$ or $1728$ and let $\widetilde{E}$ be a $k^{th}$ twist of $E$. We show that for any prime $p$ of good reduction of $\widetilde{E}$, a degree $k$ relative $p$-class group and the root number of $\widetilde{E}$ determines the dimension of the $p$-Selmer group of $\widetilde{E}$. As a consequence, we construct families of large rank $p$-class group. We also relate congruent number and cube sum problem with relative $p$-class group.
\end{abstract}

\maketitle

\section{Introduction}

The Selmer group of an elliptic curve and class group of certain number field can sometime be embedded inside a Galois cohomology group and thus one expects that Selmer group of elliptic curve and certain class group should be closely related and this relationship has been explored by various mathematicians \cite{Cas, Sch, CES, Li, sudhanshu, JMS}. In this article, we consider CM elliptic curves $E$ with $j(E) \in \{0, 1728\}$.  For the elliptic curves with $j(E)=1728$, we show that the variation of rank of $p$-Selmer group of $E$ in a quadratic (resp. quartic) twist family is completely determined by the root number and the variation of rank of certain component of $p$-class group of a family of relative quadratic (resp. quartic) extensions of a number field $K_E$ (resp. $F_E$). More precisely, let $E$ be an elliptic curve with $j(E)=1728$ and let $\widetilde{E}$ be a quadratic (resp. quartic) twist of $E$. Let $p$ be an odd prime such that both $E$ and $\widetilde{E}$ have good reduction at $p$. Let $K =\Q(i)$ and $\p$ be a prime in $K$ lying over $p$. For $C \in \{E, \widetilde{E} \}$, let
$L_{C,\p}:= K(C[\p])$ be the $\p$ division field of $C$, then $L_{C,\p}$ is Galois over $K$ with cyclic Galois group $G$ of order $N_{K/\Q}(\p)-1$. Let $K_{C,\p}$ (resp. $F_{C,\p}$) be the subfield of $L_{C,\p}$ containing $K$ such that $\Gal(L_{C,\p}/K_{C,\p}) \cong \frac{\Z}{2\Z}$ (resp. $\Gal(L_{C,\p}/F_{C,\p}) \cong \frac{\Z}{4\Z}$). Let $\chi_C$ denotes the character of $G$ which gives the action of $G$ on $C[\p]$, and let $r_\p(C)$ denotes the $k_\p:=\frac{\mathcal{O}_K}{\p}$-dimension of the $\chi_C$ component of the relative quadratic (resp. quartic) class group $\mathrm{Cl}(L_{C,\p}/K_{C,\p})[p]$ (resp. $\mathrm{Cl}(L_{C,\p}/F_{C,\p})[p]$). Finally, let $s_p(C)$ denotes the $\F_p$-dimension of $p$-Selmer group of $C$ (over $\Q$), and $\omega(C)$ denotes the root number of $C$ (over $\Q$). 

\begin{lthm}[Theorem \ref{2.3}, Theorem \ref{2.4}, Proposition \ref{2.7}]
Let $\widetilde{E}$ be a quadratic (resp. quartic) twist of an elliptic curve $E$ with $j(E)=1728$. With notation as above, $s_p(\widetilde{E})$, the rank of $p$-Selmer group of $\widetilde{E}$ is completely determined by the root number $\omega(\widetilde{E})$ and $r_\p(\widetilde{E})$, the $k_\p:=\frac{\mathcal{O}_K}{\p}$-dimension of the $\chi_{\widetilde{E}}$ component of the relative $p$-class group $\mathrm{Cl}(L_{\widetilde{E},\p}/K_{\widetilde{E},\p})[p]$ (resp. $\mathrm{Cl}(L_{\widetilde{E},\p}/F_{\widetilde{E},\p})[p]$).\\
Further, $K_{E,\p} \cong K_{\widetilde{E}, \p} $ (resp.  $F_{E,\p} \cong F_{\widetilde{E}, \p} $) and the field (up to isomorphism) $K_{E,\p}$ (resp. $F_{E,\p}$) is independent of the choice of the prime $\p$ lying above of $p$.
\end{lthm}

Similarly, let $E$ be an elliptic curve with $j(E)=0$ and let $\widetilde{E}$ be a quadratic (resp. cubic, resp. sextic) twist of $E$. Let $p$ be an odd prime such that both $E$ and $\widetilde{E}$ have good reduction at $p$. Let $K =\Q(\zeta_3)$ and $\p$ be a prime in $K$ lying over $p$. For $C \in \{E, \widetilde{E} \}$, let
$L_{C,\p}:= K(C[\p])$ be the $\p$ division field of $C$, then $L_{C,\p}$ is Galois over $K$ with cyclic Galois group $G$ of order $N_{K/\Q}(\p)-1$. Let $K^\prime_{C,\p}$ (resp. $K_{C,\p}$, resp. $F_{C,\p}$) be the subfield of $L_{C,\p}$ containing $K$ such that $\Gal(L_{C,\p}/K^\prime_{C,\p}) \cong \frac{\Z}{2\Z}$ (resp. $\Gal(L_{C,\p}/K_{C,\p}) \cong \frac{\Z}{3\Z}$, resp. $\Gal(L_{C,\p}/F_{C,\p}) \cong \frac{\Z}{6\Z}$). Let $\chi_C$ denotes the character of $G$ which gives the action of $G$ on $C[\p]$ and let $r_\p(C)$ denotes the $k_\p:=\frac{\mathcal{O}_K}{\p}$-dimension of the $\chi_C$ component of the relative quadratic (resp. cubic, resp. degree six) class group $\mathrm{Cl}(L_{C,\p}/K^\prime_{C,\p})[p]$ (resp. $\mathrm{Cl}(L_{C,\p}/K_{C,\p})[p]$, resp. $\mathrm{Cl}(L_{C,\p}/F_{C,\p})[p]$).

\begin{lthm}[Theorem \ref{3.2}, Theorem \ref{3.3}, Proposition \ref{3.6}]
Let $\widetilde{E}$ be a quadratic (resp. cubic, resp. sextic) twist of an elliptic curve $E$ with $j(E)=0$. With notation as above, the rank of $p$-Selmer group of $\widetilde{E}$ is completely determined by the root number $\omega(\widetilde{E})$ and $r_\p(\widetilde{E})$, the $k_\p:=\frac{\mathcal{O}_K}{\p}$-dimension of the $\chi_{\widetilde{E}}$ component of the relative $p$-class group $\mathrm{Cl}(L_{\widetilde{E},\p}/K^\prime_{\widetilde{E},\p})[p]$(resp. $\mathrm{Cl}(L_{\widetilde{E},\p}/K_{\widetilde{E},\p})[p]$, resp. $\mathrm{Cl}(L_{\widetilde{E},\p}/F_{\widetilde{E},\p})[p]$).\\ 
Further, $K^\prime_{E,\p} \cong K^\prime_{\widetilde{E}, \p} $ (resp. $K_{E,\p} \cong K_{\widetilde{E}, \p} $, resp.  $F_{E,\p} \cong F_{\widetilde{E}, \p} $) and the field (up to isomorphism) $K^\prime_{E,\p}$ (resp. $K_{E,\p}$, resp. $F_{E,\p}$) is independent of the choice of the prime $\p$ lying above of $p$.
\end{lthm}

In \cite{Ru1} the author showed that for an elliptic curve with $j(E) \in \{ 0, 1728 \}$, the $p$-Selmer rank of $E$ can be determined from the root number of $E$ and the structure of the $p$-class group of the $\pi$-division field of $E$. Our results can be thought of as an refinement of the result of \cite{Ru1} which is obtained by studying the structure of $\pi$-division field of $E$.\\

 In the first set of applications of our result, we give a reformulation of rank $0$ $p$-converse theorems in terms of relative $p$-class group.

\begin{lpro}[Corollary \ref{2.5}, Corollary \ref{3.4}]
    Let $E$ be an elliptic curve with $j(E) =1728$ (resp. $j(E)=0$) and $\omega(E)=+1$. The following statements are equivalent:
    \begin{enumerate}
        \item $\rk_\Q \, E =0$ and $\Sh(E/\Q)$ is finite.
        \item there exists an odd prime $p$ such that $E$ has good reduction at $p$ and for a prime $\p$ in $\Z[i]$ (resp. $\p$ in $\Z[\zeta_3]$) lying over $p$, such that $ \mathrm{Cl}(L_{E,\p}/F_{E,\p})[p](\chi_E)=0$.
        \item $L(E,1) \neq 0$.
    \end{enumerate}
\end{lpro}

Similarly, rank $1$ $p$-converse theorem for these curves can be restated as follows:

\begin{lpro}[Corollary \ref{2.6}, Corollary \ref{3.5}]
    Let $E$ be an elliptic curve with $j(E) =1728$ (resp. $j(E)=0$) and $\omega(E)=-1$. The following statements are equivalent:
    \begin{enumerate}
        \item $\rk_\Q \, E = 1$ and $\Sh(E/\Q)$ is finite.
        \item there exists an odd prime $p $ such that $E$ has good reduction at $p$, $p = \p \overline{\p}$ in $\Z[i]$ (resp. in $\Z[\zeta_3]$), such that $\mathrm{Cl}(L_{E,\p}/F_{E,\p})[p](\chi_E)$ is a cyclic group.
        \item $ord_{s=1} \, L(E,s) = 1$ .
    \end{enumerate}
\end{lpro}
As an application to the congruent number problem, we obtain
\begin{lcor} [Corollary \ref{2.8}]
    Let $n \equiv 1,2,3 \pmod8$ be a square free natural number which is a congruent number. For all odd primes $p$ co-prime to $n$, the class number of the relative quadratic extension $L_{E,\p}/K_{E,\p}$ is divisible by $p$, here $E$ is the elliptic curve $y^2=x^3-n^2x$.
\end{lcor}
Similarly, as an application to the cube sum problem, we obtain
\begin{lcor}[Corollary \ref{3.7}]
    Let $n >2$ be a cube-free integer for which the root number of $E: y^2=x^3-432n^2$ is $+1$. If there are two rational numbers $a$ and $b$ such that $n=a^3+b^3$, then for any prime $p$ co-prime to $6n$,  the class number of the relative cubic extension $L_{E,\p}/K_{E,\p}$ is divisible by $p$.
\end{lcor}


Finally, we provide some modest evidence towards the conjecture on unboundedness of $n$-rank of class groups by using families of elliptic curves with large rank.     
\begin{conj}\cite[Conjecture 1.2]{bilu2018}\label{biluconj}
  Let $n>1$ be an integer and $F$ be a number field. Then $\rk \, Cl(L)[n]$ is unbounded when  $L$ varies in the extensions of $F$ of degree $d$.  
\end{conj}

The conjecture is widely open if $d$ and $n$ are co-prime. For $2 \le d \le (p-1)$, and if $ \zeta_p \in F$, \cite[Theorem 1.7]{bilu2018} shows that there are infinitely many non-isomorphic number fields $L$ such that $[L:F]=d$ and $\rk \,Cl(L)[p] \ge 3 + \rk \, Cl(F)[p]$. Towards this, for $d=4$, we prove

\begin{lthm}[Theorem \ref{4.2}]
\begin{enumerate}
    \item Let $E:y^2=x^3 -x$, $p \equiv 3 \pmod4$ be a prime and let $F$ denote the number field $F_{E,p}$. There are infinitely many non-isomorphic number fields $L$ such that $[L:F]=4$ and $\rk \,Cl(L)[p] \ge 12 + \rk \, Cl(F)[p]$. 
    \item Let $E:y^2=x^3 -x$, $p \equiv 1 \pmod4$ be a prime, $\p$ a prime in $\Z[i]$ lying above $p$ and let $F$ denote the number field $F_{E,\p}$. There are infinitely many non-isomorphic number fields $L$ such that $[L:F]=4$ and $\rk \,Cl(L)[p] \ge 6 + \rk \, Cl(F)[p]$. 
\end{enumerate}   
\end{lthm}

We remark that if $\dim_{\F_p} \Sel^p(E/\Q)$ is unbounded as $E$ varies over the CM elliptic curves with $j(E)=1728$, we'll obtain that the Conjecture \ref{biluconj} is true for $d=4$ and $F \cong F_{E,\p}$.
Similarly, for $d=6$, we show
\begin{lthm}[Theorem \ref{4.5}]
\begin{enumerate}
    \item Let $E:y^2=x^3 +1$, $p \equiv 2 \pmod 3$ be a prime and let $F$ denote the number field $F_{E,p}$. There are infinitely many non-isomorphic number fields $L$ such that $[L:F]=6$ and $\rk \,Cl(L)[p] \ge 12 + \rk \, Cl(F)[p]$. 
    \item Let $E:y^2=x^3+1$, $p \equiv 1 \pmod3$ be a prime, $\p$ a prime in $\Z[\zeta_3]$ lying above $p$ and let $F$ denote the number field $F_{E,\p}$. There are infinitely many non-isomorphic number fields $L$ such that $[L:F]=6$ and $\rk \,Cl(L)[p] \ge 6 + \rk \, Cl(F)[p]$. 
\end{enumerate}   
\end{lthm}

Again, we remark that if $\dim_{\F_p} \Sel^p(E/\Q)$ is unbounded as $E$ varies over the CM elliptic curves with $j(E)=0$, we'll obtain that the Conjecture \ref{biluconj} is true for $d=6$ and $F \cong F_{E,\p}$. \\

The paper is structured as follows: In Section \ref{S2}, we study the $p$-Selmer group of elliptic curves with $j$-invariant $1728$. In Section \ref{S3}, we study the $p$-Selmer group of elliptic curves with $j$-invariant $0$. In Section \ref{S4}, we construct class group with large $p$-rank. In the final Section \ref{nc}, we use SAGE to compute class groups and use that information to determine the rank of various elliptic curves with $j$-invariant $0$ and $1728$.\\

\noindent {\bf Acknowledgements:}  We thank Ashay Burungale and B. Sury for discussions.

\section{$p$-Selmer group of $E_D: Y^2=X^3-DX$}\label{S2}
Let $E:=E_1$ denote the elliptic curve $E: Y^2=X^3-X$. For a nonzero integer $D \in \Z$, the elliptic curve $E_D: Y^2=X^3-DX$ is the quartic twist of $E$ by $D$. The elliptic curve $E_D$ has CM by $\Z[i]$ given as $i(x,y)=(-x, iy)$ and for a prime $q \nmid 2D$, $E_D$ has good ordinary reduction at $q$ if and only if $q \equiv 1 \pmod{4}$. By $\omega(E_D)$, we denote the global root number of $E_D$.\\

The following lemmas provide an explicit description of $n$-division polynomial of $E_D$ which will be useful for our study.

\begin{lemma}\label{2.1}
    For every integer $n \ge 1$, there exists a homogeneous polynomial $f_n(Z,W) = c_n Z^{s_n}+ \cdots + d_n W^{s_n} \in \Z[Z,W]$, where $ (c_n, d_n, s_n)  = \begin{cases}
        (n, (-1)^{\frac{n-1}{2}} , \frac{n^2-1}{4}) & \text{ if } n \text{ is odd,}\\
        (\frac{n}{2}, \frac{n}{2}, \frac{n^2-4}{4})   & \text{ if } n \text{ is even}, \end{cases} $ 
        
       \noindent such that
    $$ \psi_n= \begin{cases}
        f_n(x^2,D) & \text{ if } n \text{ is odd,}\\
        2yf_n(x^2, D) & \text{ if } n \text{ is even}.
    \end{cases} $$
    Here $\psi_n \in \Z[x,y]$ is the $n^{th}$ division polynomial for $E_D: y^2=x^3-Dx$.  
\end{lemma}

\begin{proof}
Note that $f_1(Z,W)=f_2(Z,W)=1$, $f_3(Z,W)=3Z^2-6ZW-W^2$, $f_4(Z,W)=2(Z^3-5WZ^2-5W^2Z+W^3)$. We prove the result by induction on $n$ and using the recursive relations that define $\psi_n$.\\
Assume that the statement is true for $1 \le t \le n-1$. Let $n = 4m+1$, then $\psi_n=\psi_{2m+2}\psi_{2m}^3-\psi_{2m-1}\psi_{2m+1}^3$, which implies that $\psi_n=16y^4f_{2m+2}(x^2,D)f_{2m}^3(x^2,D)-f_{2m-1}(x^2,D)f_{2m+1}^3(x^2,D)$. We define $f_{4m+1}(Z,W)=16Z(Z-W)^2f_{2m+2}(Z,W)f_{2m}^3(Z,W)- f_{2m-1}(Z,W)f_{2m+1}^3(Z,W)$ and using the induction on degree of $f_t(Z,W)$, one checks that $f_{4m+1}(Z,W)$ is a homogeneous polynomial in $Z$ and $W$ of degree $ \frac{(4m+1)^2-1}{4}$ and we have $\psi_{4m+1} = f_{4m+1}(x^2,D)$. Similarly, one can show that $f_{4m+3}(Z,W) =f_{2m+3}(Z,W)f_{2m+1}^3(Z,W) - 16Z(Z-W)^2f_{2m}(Z,W)f_{2m+2}^3(Z,W) $ and $f_{2m}(Z,W)= f_m(Z,W)[f_{m+2}(Z,W) f_{m-1}^2(Z,W) - f_{m-2}(Z,W)f_{m+1}^2(Z,W)]$. Also the formulas of $s_n$ and $d_n$ can be proved using induction. Finally, since $$\psi_n = \begin{cases}
      nx^{\frac{n^2-1}{2}} + \cdots, & \text{ if } n \text{ is odd},\\
      y(nx^{\frac{n^2-4}{2}}+ \cdots), & \text{ if } n \text{ is even},
    \end{cases}$$
    it follows that $c_n=n$ if $n$ is odd and $c_n=\frac{n}{2}$ if $n$ is even. \qedhere
\end{proof}

\begin{lemma}\label{2.2}
    For every integer $n \ge 1$, there exists a homogeneous polynomial $g_n(Z,W) = Z^{s_n}+ \cdots \in \Z[Z,W]$, where  $s_n = \begin{cases}
        \frac{n^2-1}{4} & \text{ if } n \text{ is odd,}\\
        \frac{n^2}{4} & \text{ if } n \text{ is even},
    \end{cases}$ 
    
    \noindent such that
    $$ 
    \phi_n= \begin{cases}
        xg_n(x^2,D)^2 & \text{ if } n \text{ is odd,}\\
        g_n(x^2, D)^2 & \text{ if } n \text{ is even}.
    \end{cases} $$
    Here $\phi_n= x\psi_n^2 - \psi_{n-1}\psi_{n+1}  \in \Z[x]$ and $\psi_n$ is the $n^{th}$ division polynomial for $E_D: y^2=x^3-Dx$.  
\end{lemma}

\begin{proof}
Using induction as in the proof of the preceding lemma, one can show that there exists homogeneous polynomial $\tilde{g}_n(Z,W)$ of degree $2s_n$ such that
$$ 
    \phi_n= \begin{cases}
        x\tilde{g}_n(x^2,D) & \text{ if } n \text{ is odd,}\\
        \tilde{g}_n(x^2, D) & \text{ if } n \text{ is even.}
    \end{cases} $$
Now from the point addition formula of elliptic curve, observe that $\frac{\phi_{2n}(P)}{\psi_{2n}(P)^2} = x(2nP) = \frac{\phi_2(nP)}{\psi_2(nP)^2}$ and hence $\phi_{2n}(P) = \phi_2(nP) \Big(\frac{\psi_{2n}(P)}{\psi_2(nP)}\Big)^2$. For the elliptic curve $E_D$, $\phi_2(P) = (x(P)^2+D)^2$ and hence 
    $$ \phi_{2n}(P) = \Bigg(\frac{(x(nP)^2+D)\psi_{2n}(P) }{\psi_2(nP)} \Bigg)^2 = [ \phi_n^2(P) + D \psi_n^4(P)]^2,$$
    equivalently, $\tilde{g}_{2n}(x^2,D) = g_{2n}(x^2,D)^2 \in \Z[x]$.\\
On the other hand, if $n$ is odd, then from the formula $4y\omega_n = \psi_{n+2} \psi_{n-1}^2 - \psi_{n-2}\psi_{n+1}^2$, it follows that $\omega_n = y \tilde{f}_n(x^2,D)$, where $\tilde{f}_n(Z,W) \in \Z[Z,W]$ is a homogeneous polynomial. Now consider the polynomial $h_n(Z,W):= Z\tilde{g}_n(Z,W)^2 -Wf_n(Z,W)^4 \in \Z[Z,W]$, then $\tilde{g}_n(x^2,D) h_n(x^2,D) = (x^2-D)\tilde{f}_n(x^2,D)^2$ as we have $\omega_n^2= \phi_n^3 - D \phi_n \psi_n^4 = \phi_n(\phi_n^2 - D \psi_n^4)$. Equivalently, in $\Z[X]$, we have the polynomial identity $\tilde{g}_n(X^2,D) h_n(X^2,D) = (X^2-D)\tilde{f}_n(X^2,D)^2$. Since $\phi_a$ and $\psi_a$ are co-prime, it follows that  $\tilde{g}_n(X^2,D)$ and $h_n(X^2,D)$ have no common root in $\QQ$. Now suppose that $(X \pm \sqrt{D})$ is a factor of $\tilde{g}_n(X^2,D) \in \QQ[X]$, but $(X  \mp \sqrt{D})$ is not a factor of $\tilde{g}_n(X^2,D) \in \QQ[X]$. Since, for the polynomial $(x^2-D)\tilde{f}_n(x^2,D)^2 \in \QQ[X]$, multiplicity of the roots $\pm \sqrt{D}$ are odd and for all other roots, the multiplicity is even, it follows that in the above mentioned scenario, $\deg(\tilde{g}_n(X^2,D))$ is odd, which is a contradiction. Now suppose that $(X^2-D)$ is a factor of $\tilde{g}_n(X^2,D)$, and write $\tilde{g}_n(X^2,D) = (X^2-D)p_n(X^2) \in \Z[X]$. Then it follows that $p_n(X^2) = q_n(X^2)^2$ where $q_n(X) \mid \tilde{f}_n(X,D) \in \Z[X]$, which in turn implies that $\deg(\tilde{g}_n(X^2,D)) \equiv 2 \pmod{4}$, contradiction. Thus, $\tilde{g}_n(X^2,D) = g_n(X^2)^2 \in \Z[X]$, where $g_n(X) \mid \tilde{f}_n(X,D) \in \Z[X]$. \qedhere

\end{proof}

Let us first fix a prime $p \equiv 3 \pmod{4}$, and for an integer $D$ which is co-prime to $p$, let $s_p(E_D) := \dim_{\F_p} \Sel^p(E_D/\Q)$. Let $K=\Q(i)$ and $L_{D,p}:= K(E_D[p])$ be the $p$-division field of $E_D$ over $K$, then we have $G:=\Gal(L_{D,p}/K) \cong \frac{\Z}{(p^2-1)\Z}$. The following theorem shows that the variation of $p$-Selmer rank of quartic twist by $D$ of $E_1$ is entirely controlled by the relative $p$-class group of a cyclic quartic extension of a fixed number field $F_{1,p}$ and root number of $E_D$.

 \begin{theorem}\label{2.3}
  Let $p \equiv 3 \pmod{4}$ be a prime and $D$ be an integer co-prime to $p$. Let $K=\Q(i)$, $L_{D,p}:= K(E_D[p])$ be the $p$-division field of $E_D$ and $F_{D,p}$ be the field $K \subset F_{D,p} \subset L_{D,p}$ with $G=\Gal(L_{D,p}/F_{D,p}) \cong \frac{\Z}{4\Z}$. Let $\chi_D$ denotes the character which gives the action of $\Gal(L_{D,p}/K)$ on $E_D[p] \cong \F_{p^2}$ and $r_p(E_D)= \dim_{\F_{p^2}} Cl(L_{D,p}/F_{D,p})[p](\chi_D)$. Then we have,
   $$r_p(E_D) \le s_p(E_D) \le 1 + r_p(E_D),$$ 
   and hence $\omega(E_D)$ and $r_p(E_D)$ uniquely determines $s_p(E_D)$. Moreover, we have $F_{D,p} \cong F_{1,p} \cong \frac{K[X]}{(f_p(X,1))}$ and $L_{D,p} \cong \frac{K[X]}{(f_p(X^4,D))}$.
 \end{theorem}

\begin{proof}
From \cite[Theorem 1]{Ru1}, we see that $r_p(E_D) \le s_p(E_D) \le 1 + r_p(E_D)$, where $r_p(E_D)= \dim_{\F_{p^2}}\Hom(Cl(L_{D,p}),E_D[p])^{G} = \dim_{\F_{p^2}} Cl(L_{D,p})[p](\chi_D)$. Now as $p \nmid [L_{D,p}:K]$, we have $Cl(L_{D,p})[p] \cong Cl(L_{D,p}/F_{D,p})[p] \times Cl(F_{D,p})[p]$ and as $\Gal(L_{D,p}/F_{D,p})$ acts trivially on $Cl(F_{D,p})[p]$, it follows that $r_p(E_D) = \dim_{\F_{p^2}} Cl(L_{D,p}/F_{D,p})[p](\chi_D)$. From $p$-parity conjecture, we know that $s_p(E_D) \equiv \omega(E_D) \pmod{2}$, and hence $\omega(E_D)$ and $r_p(E_D)$ uniquely determines $s_p(E_D)$.\\
Let $P_D:=(x(P_D), y(P_D))$ be a non-trivial element of $E_D[p]$, then $x(P_D)$ is a root of $\psi_p = f_p(X^2,D)$ (as in Lemma \ref{2.1}). Hence $x(P_D)^2$ is a root of $f_p(X,D)$ which is a polynomial of degree $\frac{p^2-1}{4}$ and, hence $[K(x(P_D)^2):K] \le \frac{p^2-1}{4}$. Since $[K(x(P_D)) : K(x(P_D)^2)] \le 2$, $[K(y(P_D)) : K(x(P_D))] \le 2$  and $[K(E_D[p]) : K] = p^2-1$, it follows that all the inequalities must be equality and, hence $f_p(X,D)$ is irreducible in $K[X]$ and,
$K(x(P_D)^2) \cong \frac{K[X]}{(f_p(X,D))} := F_{D,p}$. Moreover, as $f_p(X,D)$ is a homogeneous polynomial, it follows that for every $D \ge 1$, $\frac{x(P_D)^2}{D}$ is a root of $f_p(X,1)$ and, hence $F_{D,p} := K(x(P_D)^2) = K(\frac{x(P_D)^2}{D}) \cong K(x(P_1)^2):=F_{1,p}$ is independent of $D$.\\
Next, observe that as $[K(x(P_D)): K] = \frac{p^2-1}{2}$. Since $x(P_D)$ is a root of $f_p(X^2,D)$, it follows that $f_p(X^2,D)$ is irreducible in $K[X]$ and, we have $K(x(P_D))\cong \frac{\Q(i)[X]}{(f_p(X^2,D))} :=K_{D,p}$. Further note that $L_{D,p}$ is a quadratic extension of $K_{D,p}$.\\
Finally, as $(p-1)P_D = - P_D$, it follows that $ \frac{\phi_{p-1}(P_D)}{\psi_{p-1}(P_D)^2} = x((p-1)P_D) = x(-P_D) = x(P_D)$ or equivalently, 
    $$ \frac{g_{p-1}(x(P_D)^2, D)^2 }{4 y(P_D)^2 f_{p-1}(x(P_D)^2, D)^2 } = x(P_D).$$
    As a consequence we see that $L_{D,p}:= K(E_D[p]) \cong K\left(\sqrt{x(P_D)}\right) \cong \frac{K[X]}{(f_p(X^4,D))}$. \qedhere
\end{proof}

 Next we consider the case $p \equiv 1 \pmod 4$ and let $p \OO_K = \p \overline{\p}$, where $\p$ and $\overline{\p}$ are the primes in $\OO_K$ lying above $p$. Let $\pi= a+ib$ be the element in $\Z[i]$ such that $\pi \equiv 1 \pmod{2+2i}$ and, $\p = (\pi)$ and, for every natural number $D$ co-prime to $p$, we define $\pi_D = \overline{\big( \frac{D}{\pi}\big)_4} \pi$, where $\big( \frac{\cdot}{\cdot} \big)_4$ denotes the $4^{th}$ power residue symbol. With this notation, we have $\psi_D(\p) = \pi_D$ (resp. $\psi_D(\overline{\p}) = \overline{\pi_D}$), where $\psi_D$ is the Gr\"ossencharacter associated to $E_D$. From now on, fix a prime $\p$ in $\OO_K$ lying above $p$.   Let $L_{D, \p} := K(E_D[\pi_D])$ and
$G:= \Gal(L_{D,\p}/K) \cong \frac{\Z}{(p-1)\Z}$. Finally, define $f_{D,p}(X)= 4X(X-D)g_{|a|}(X,D)^2f_{|b|}(X,D)^2+g_{|b|}(X,D)^2f_{|a|}(X,D)^2 = pX^{\frac{p-1}{2}}+ \cdots + D^{\frac{p-1}{2}}$, where $\pi = a+ib \equiv 1 \pmod{2+2i}$, and the polynomials $f_n,g_n$ are as in Lemma \ref{2.1},\ref{2.2} respectively.  The following theorem shows that the variation of $p$-Selmer rank of quartic twist by $D$ of $E_1$ is entirely controlled by the relative $p$-class group of a cyclic quartic extension of a fixed number field $F_{1,p}$ and the root number of $E_D$.

\begin{theorem}\label{2.4}
    Let $p \equiv 1 \pmod{4}$ be a prime with $(p,D)=1$. Let $K=\Q(i)$ and $\p$ be a prime in $\OO_K$ lying above $p$. Let $L_{D, \p}:= K(E_D[\pi_D])$ and $F_{D, \p}$ be the unique subfield of $L_{D,\p}$ which contains $K$ and satisfies $\Gal(L_{D,\p}/F_{D,\p}) \cong \frac{\Z}{4\Z}$.  Let $\chi_D$ be the character which gives the action of $G$ on $E_D[\pi_D] \cong \F_p$ and $r_\p(E_D):= \dim_{\F_p} Cl(L_{D,\p}/F_{D, \p})[p](\chi_D)$. Then we have
    $$r_\p(E_D) \le s_p(E_D) := \dim_{F_p} {\mathrm{Sel}}^p(E_D/\Q) \le 1+ r_\p(E_D),$$
    and hence $\omega(E_D)$ and $r_\p(E_D)$ uniquely determines $s_p(E_D)$. Further, we have $F_{D, \p} \cong \frac{\Q[X]}{(f_{1,p}(X))}$ and $L_{D, \p} \cong \frac{\Q[X]}{(f_{D,p}(X^4))}$.
\end{theorem}

\begin{proof}
    From \cite[Theorem 1]{Ru1}, we see that $r_\p(E_D) \le s_p(E_D) \le 1 + r_\p(E_D)$, where $r_\p(E_D)= \dim_{\F_{p}}\Hom(Cl(L_{D,\p}),E_D[\pi_D])^{G} = \dim_{\F_{p}} Cl(L_{D,\p})[p](\chi_D)$. Now as $p \nmid [L_{D,\p}:K]$, we have $Cl(L_{D,\p})[p] \cong Cl(L_{D,\p}/F_{D,\p})[p] \times Cl(F_{D,\p})[p]$ and, since $\Gal(L_{D,\p}/F_{D,\p})$ acts trivially on $Cl(F_{D,\p})[p]$, it follows that $r_\p(E_D) = \dim_{\F_{p}} Cl(L_{D,\p}/F_{D,\p})[p](\chi_D)$. From $p$-parity conjecture, we know that $s_p(E_D) \equiv \omega(E_D) \pmod{2}$, and hence $\omega(E_D)$ and $r_\p(E_D)$ uniquely determine $s_p(E_D)$.\\
    Now let $P_D=(x(P_D), y(P_D))$ be a generator of $E_D[\pi_D]$. Note that $\pi_D = a_D + ib_D = i^{k}(a+ib)$ and hence $(|a_D|, |b_D|)$ is either $(|a|, |b|)$ or $(|b|, |a|)$. Since $\pi_DP_D=0$, it follows that $x(|a_D|P_D) = - x(|b_D|P_D)$ and hence from Lemma \ref{2.1} and Lemma \ref{2.2}, it follows that
    $$ \frac{x(P_D)g_{|a|}(x(P_D)^2,D)^2}{ f_{|a|}(x(P_D)^2,D)^2} = - \frac{g_{|b|}(x(P_D)^2,D)^2}{4 y(P_D)^2 f_{|b|}(x(P_D)^2,D)^2},$$
    and hence $K(E_D[\pi_D]) = K(x(P_D), y(P_D)) \cong K(\sqrt{x(P_D)})$. Again from $\pi_DP_D=0$, it follows that $x(P_D)^2$ is a root of $f_{D,p}(X)$, which is a homogeneous polynomial in $D$ and $X$ and hence it follows that $K(x(P_D)^2)$ is independent of $D$ (upto isomorphism). Since $[K(\sqrt{x(P_D)}) : K(x(P_D)^2)] \le 4$, it follows that $F_{D, \pi_D} \subset K(x(P_D)^2)$. We will show that $F_{D, \pi_D} = K(x(P_D)^2)$. Note that  $f_{D,p}(X^4)$ can not be irreducible over $K$, as otherwise, $[K(\sqrt{x(P_D)}):K] = 2(p-1)$ which is a contradiction. Consequently, $f_{D,p}(X^4)$ has at-least one irreducible factor of degree $(p-1)$ over $K$.\\
    By \cite[Lemma 4]{CW}, we know that the conductor of $L_{D,\p}$ over $K$ is $\p f_D$, where $f_D$ is the conductor of Gr\"o{\ss}encharacter $\psi_D$ of $E_D$, in-particular, $\overline{\p}$ is unramified in $L_{D,\p}/K$. Further by \cite[Lemma 5]{CW}, we also know that $\p$ is totally ramified in $L_{D,\p}/K$.\\
    Similarly, writing $Q_D=(x(Q_D), y(Q_D))$ as a generator of $E_D[\overline{\pi_D}]$, we see that $L_{D,\overline{\p}} := K(E_D[\overline{\pi_D}]) \cong K(\sqrt{x(Q_D)})$, $\p$ is unramified in $L_{D, \overline{\p}}/K$ and $\overline{\p}$ is totally ramified in $L_{D, \overline{\p}}/K$. In-particular, we see that $L_{D,\p} \cap L_{D, \overline{\p}} =K$. Further, as $L_{D,\p} \cong L_{D, \overline{\p}}$, we see that for any isomorphism $\theta: L_{D,\p} \to L_{D, \overline{\p}}$, we have $\theta|_K$ is the non-trivial automorphism of $K$. Now if $f_{D,p}(X^4)$ has a unique irreducible factor, say $h(X)$, of degree $(p-1)$ over $K$, then both $\sqrt{x(P_D)}$ and $\sqrt{x(Q_D)}$ are roots of $h(X)$ and there exists an isomorphism $\theta:L_{D,\p} \to L_{D, \overline{\p}}$ with $\theta(\sqrt{x(P_D)}) = \sqrt{x(Q_D)}$. Consequently, $h(X)=\overline{h(X)}$, which contradicts the fact $\theta|_K$ is the complex conjugation. Hence, $f_{D,p}(X^4) = h_1(X)h_2(X) \in K[X]$, where both $h_1(X)$ and $h_2(X)$ are irreducible polynomials of degree $p-1$ and $\sqrt{x(P_D)}$ is a root of $h_1(X)$ (say) and $\sqrt{x(Q_D)}$ is a root of $h_2(X)$, and hence $L_{D,\p} \cong \frac{K[X]}{(h_1(X))}$ and $L_{D,\overline{\p}} \cong \frac{K[X]}{(h_2(X))}$. Since $L_{D,\p} \cong L_{D, \overline{\p}}$ which induces complex conjugation when restricted to $K$, it follows that $h_2(X) = \overline{h_1(X)}$ and consequently, $f_{D,p}(X^4) = h_1(X) \overline{h_1(X)} \in K[X]$, which in-turn implies that $f_{D,p}(X^4)$ is irreducible of degree $2(p-1)$ over $\Q$ and $L_{D,\p} \cong \frac{\Q[X]}{(f_{D,p}(X^4))}$. Since, $f_{D,p}(X^4)$ is irreducible over $\Q$, it follows that $f_{D,p}(X)$ is irreducible over $\Q$. Again by degree argument we can say that $f_{D,p}(X)$ can not be irreducible over $K$. Using the above argument, it follows that we have an isomorphism $K(x(P_D)^2) \cong K(x(Q_D)^2)$ which induces complex conjugation on $K$, which in-turn implies that $f_{D,p}(X) = h_3(X)\overline{h_3(X)} \in K[X]$ and consequently, $F_{D,\p} \cong \frac{\Q[X]}{(f_{D,p}(X))} \cong \frac{\Q[X]}{(f_{1,p}(X))}$. \qedhere
    
\end{proof}

As a consequence of Theorem \ref{2.3} and Theorem \ref{2.4}, we obtain

\begin{corollary}\label{2.5}
  Suppose that $\omega(E_D) = +1$. Then, $L(E_D, 1) \neq 0$ $\iff$ $\rk_\Q E_D =0$ and $\Sh(E_D/\Q)$ is finite   $\iff$ for all but finitely many primes $\p$ in $K$ with $\p \nmid 2D$, $r_\p(E_D) =0$ $\iff$ there exists a prime $\p$ in $K$ with $\p \nmid 2D$ with $r_\p(E_D) =0$ .
\end{corollary}

\begin{proof}
    First suppose that $L(E_D,1) \neq 0$, then by the work of Coates-Wiles \cite{CW} and Rubin \cite{Ru3}, it follows that $\rk_\Q E_D =0$ and $\Sh(E_D/\Q)$ is finite, which in-turn implies that for all but finitely many $p$, we have $\Sel^p(E_D/\Q) = 0$. Consequently by either Theorem \ref{2.3} or Theorem \ref{2.4}, for all but finitely many $\p$, we have $r_\p(E_D)=0$. Now if $r_\p(E_D)=0$, then by either Theorem \ref{2.3} or Theorem \ref{2.4}, we have $s_p(E_D)=0$ and, by $p$-converse theorem of Rubin \cite{Ru2}, it follows that $L(E_D,1) \neq 0$. \qedhere
\end{proof}

In a similar fashion, we prove:
\begin{corollary} \label{2.6}
    Suppose that $\omega(E_D) = -1$. Then, $L(E_D, 1) = 0$ and $L^\prime(E_D,1) \neq 0$  $\iff$ $\rk_\Q E_D =1$ and $\Sh(E_D/\Q)$ is finite    $\iff$  for all but finitely many primes $\p$ in $K$ with $\p \nmid 2D$ with $r_\p(E_D) \le 1$ $\iff$ there exists a prime $\p \nmid 2D$ in $K$ lying over $p \equiv 1 \pmod 4$  with $r_\p(E_D) \le 1$.  
\end{corollary}

\begin{proof}
   First suppose that $L(E_D,s)$ has simple zero at $s=1$, then by the work of Gross-Zagier \cite{GZ} and Kolyvagin \cite{Kol}, it follows that $\rk_\Q E_D =1$ and $\Sh(E_D/\Q)$ is finite, which in-turn implies that for all but finitely many $p$, we have $\dim_{F_p} \Sel^p(E_D/\Q) = 1$. Consequently by either Theorem \ref{2.3} or Theorem \ref{2.4}, for all but finitely many $\p$, we have $r_\p(E_D) \le 1$. Now if for a prime $\p \nmid 2D$ in $K$ lying over $p \equiv 1 \pmod 4$ we have $r_\p(E_D) \le 1$, then by Theorem \ref{2.4}, we have $s_p(E_D)=1$ and by $p$-converse theorem of Burungale-Tian \cite{BT}, it follows that $\mathrm{ord}_{s=1}  L(E_D,s)=1$. \qedhere
\end{proof}

Observe that for every $D$ and a prime $\p \nmid 2D$  in $K$, we have a chain of sub-fields $K \subset F_{1,\p}  \subsetneq K_{D,\p} \subsetneq  L_{D,\p}$. Note that we can explicitly describe $K_{D,\p}$ as $K_{D,p} \cong \frac{K[X]}{(f_p(X^2,D))}$ if $p \equiv 3 \pmod 4$, and $K_{D,\p} \cong \frac{\Q[X]}{(f_{D,p}(X^2))}$ if $\p \mid p$ with $p \equiv 1 \pmod 4$. Further note that if $D_1 = D\alpha^4$ for some $\alpha \in K$, $L_{D_1,\p} = L_{D,\p}$ and if $D_2 = D \beta^2$ for some $\beta \in K$, then $K_{D_2,\p} = K_{D,\p}$. Consequently the following result shows that for quadratic twists of $E_D$, the variation of $p$-Selmer rank is controlled by a relative quadratic extension and the root number.

\begin{proposition}\label{2.7}
   Let $\p \nmid 2D$ be a prime in $K$ lying over $p$. Let $M$ be an integer co-prime to $pD$ and $E_D^M: Y^2=X^3-DM^2X$ denote the quadratic twist of $E_D$ by $M$. Then we have,
    $r_\p(E_D^M) \le s_p(E_D^M) \le 1 + r_\p(E_D^M)$, where $r_\p(E_D^M)= \dim_{\mathcal{O}_K/ \p} Cl(L_{DM^2,\p}/K_{DM^2,\p})[p](\chi_{DM^2})$.  
 Moreover, we have $K_{DM^2,\p} \cong K_{D,\p}$.  
\end{proposition}

In-particular, for congruent number elliptic curves we obtain:

\begin{corollary}\label{2.8}
    Let $n \equiv 1,2,3 \pmod{8}$ be a square-free integer. Suppose there exists a prime $\p \nmid 2n$ such that $r_\p(E_1^n) =0$, then $n$ is not a congruent number.\\
    In particular, if a square-free integer $n \equiv 1,2,3 \pmod{8}$ is a congruent number, then for all primes $\p \nmid 2n$, we have $p \mid |Cl(L_{n^2,\p}/K_{1,\p})|$.
\end{corollary}

\begin{proof}
    Since $n \equiv 1,2,3 \pmod{8}$, it follows that $\omega(E_1^n) =+1$ and, hence $r_\p(E_1^n) =0$ implies that $s_p(E_1^n) =0$, which in-turn implies that $\rk_\Q E_1^n =0$, that is $n$ is not a congruent number.\\
    Note that, $n \equiv 1,2,3 \pmod{8}$ is a congruent number which implies that $s_p(E_1^n) \ge 2$ or, equivalently $r_p(E_1^n) \ge 1$, which in-turn implies that $p \mid |Cl(L_{n^2,\p}/K_{1,\p})|$. \qedhere
\end{proof}

Given a non-negative integer $k$, by $D(k)$ we denote the integer $D(k)= (8k+5)(8k+6)(8k+7)$ and by $D(k)_{sf}$ we denote the square-free part of $D(k)$. One checks that for every $k \ge 0$, $D(k)_{sf}$ is a congruent number and, $D(k)_{sf} \equiv 2 \pmod 8$. Now for a prime $p >5$ and $p \equiv 1 \pmod 4$, one obtains that the infinite set $S_p=\{ D(pk)_{sf}\}_{k \ge 0}$ consists of congruent numbers which are co-prime to $p$. Similarly the infinite set $S_5 = \{ D(5k+2)_{sf}, D(5k+4)_{sf}\}_{k \ge 0}$ consists of congruent numbers which are co-prime to $5$.\\
Next observe that, for every integer $s$, the integer  $f(s)= (2s+1)^4+24(2s+1)^2+16$ is a congruent number, and congruent to $1$ modulo $8$. For a prime $p \equiv 3 \pmod 4$, the infinite set $T_p = \{ f(ps)_{sf}\}_{s \ge 0}$ consists of congruent numbers which are co-prime to $p$. Consequently, we obtain the following:

\begin{corollary}\label{2.9}
    Given any prime $p \ge 3$, there are infinitely many quadratic extensions $L_\p$ of $K_{1,\p}:=K_\p$ such that $p$ divides the relative class number $L_\p/K_\p$.
\end{corollary}

In a similar fashion to Corollary \ref{2.8}, we obtain:
\begin{corollary}\label{2.10}
Let $n \equiv 5,6,7 \pmod{8}$ be a square-free integer.  Suppose that there exists a prime $p \equiv 1 \pmod{4}$ which is co-prime to $n$ such that $r_\p(E_1^n) \le 1$, then $n$ is a congruent number.
\end{corollary}




\section{$p$-Selmer group of $E_A: Y^2=X^3+A$}\label{S3}
Let $E:=E_1$ denote the elliptic curve $E: Y^2=X^3+1$. For a nonzero integer $A \in \Q^\times/(\Q^{\times})^6$, the elliptic curve $E_A: Y^2=X^3+A$ is the sextic twist of $E$ by $A$. The elliptic curve $E_A$ has CM by $\Z[\zeta_3]$, given as $\zeta_3(x,y) = (\zeta_3 x, y)$, and $E_A$ has good reduction at $p$ for $p \nmid 6A$. Hence by Deuring's theorem, for primes $p \nmid 6A$,  $E_A$ has good ordinary reduction if and only if $p \equiv 1 \pmod{3}$. By $\omega(E_A)$, we denote the global root number of $E_A$. 

Following lemma provides description of the $n$-th division polynomial for the elliptic curve $E_A$.
\begin{lemma}\label{cubesumpsin}
    For every integer $n \ge 1$, there exists a homogeneous polynomial $f_n(Z,W) = c_n Z^{d_n}+ \cdots \in \Z[Z,W]$, where $ c_n = \begin{cases}
        n & \text{ if } n \text{ is odd,}\\
        \frac{n}{2} & \text{ if } n \text{ is even},
    \end{cases}$ and the degree $d_n$ is given by
    $$ d_n= \begin{cases}
    \frac{n^2 -1}{6} & \text{ if } n \equiv 1, 5 \pmod{6},\\
    \frac{n^2-3}{6} & \text{ if } n \equiv 3 \pmod{6},\\
    \frac{n^2-4}{6}  & \text{ if } n \equiv 2, 4 \pmod{6},\\
    \frac{n^2-6}{6}  & \text{ if } n \equiv 0 \pmod{6},    
    \end{cases}
$$ 
such that
    $$ 
    \psi_n= \begin{cases}
        f_n(y^2,A) = f_n(x^3+A, A) & \text{ if } n \equiv 1, 5 \pmod{6},\\
        xf_n(y^2, A) = xf_n(x^3+A, A) & \text{ if } n \equiv 3 \pmod{6},\\
        2yf_n(y^2,A)= 2yf_n(x^3+A,A) & \text{ if } n \equiv 2, 4 \pmod{6},\\
        2xyf_n(y^2,A)= 2xyf_n(x^3+A,A) & \text{ if } n \equiv 0 \pmod{6}.
    \end{cases} $$
    Here $\psi_n$ is the $n^{th}$ division polynomial for $E_A: y^2=x^3+A$.  
\end{lemma}

\begin{proof}
    Note $f_1(Z,W)= f_2(Z,W)=1$, $f_3(Z,W)= 3(Z+3W)$, $f_4(Z,W)= 2(Z^2+18ZW-27W^2)$. One can prove the result by induction on $n$ and using the recursive relations that define $\psi_n$. \qedhere
\end{proof}


Let $p \equiv 5 \pmod{6}$ be a prime, and for an integer $A$ co-prime to $p$, define $s_p(E_A) := \dim_{\F_p} \Sel^p(E_A/\Q)$. Let $K=\Q(\zeta_3)$ and $L_{A,p}:= K(E_A[p])$ be the $p$-division field of $E_A$ over $K$, which is of degree $12d_p= 2(p^2-1)$. Also we have $G:=\Gal(L_{A,p}/K) \cong \frac{\Z}{(p^2-1)\Z}$. 

\begin{theorem}\label{3.2}
    Let $p \equiv 5 \pmod{6}$ be a prime and $A$ be an integer such that $(p,A)=1$. Let $F_{A,p}$ be the subfield such that $K:= \Q(\zeta_3) \subset F_{A,p} \subset L_{A,p} := K(E_A[p])$ with $\Gal(L_{A,p}/F_{A,p}) \cong \Z/{6\Z}$.  Let $\chi_A$ denote the character that gives the action of $G =\Gal(L_{A,p}/K)$ on $E_A[p] \cong \F_{p^2}$ and $r_p(E_A)= \dim_{\F_{p^2}} Cl(L_{A,p}/F_{A,p})[p](\chi_A)$. Then we have,
   $$r_p(E_A) \le s_p(E_A)\le 1 + r_p(E_A),$$ 
   and hence $\omega(E_A)$ and $r_p(E_A)$ uniquely determine $s_p(E_A)$.

   Furthermore, let $K_{A,\p}$ (resp. $K^\prime_{A,\p}$)  be the subfield of $L_{A,\p}$ such that $ F_{A,\p} \subset K_{A,\p} \subset L_{A,\p}$ (resp. $ F_{A,\p} \subset K^\prime_{A,\p} \subset L_{A,\p}$)  with $\Gal(L_{A,\p}/K_{A,\p}) \cong \Z/{3\Z}$ (resp. $\Gal(L_{A,\p}/K^\prime_{A,\p}) \cong \Z/{2\Z}$), then
  $$ K_{A,p} \cong \frac{K[Y]}{(f_p(Y^2,A))}, K^\prime_{A,p} \cong \frac{K[X]}{(f_p(X^3+A,A))} \text{ and }  F_{A,p} = K_{A,p} \cap K^\prime_{A,p} \cong \frac{K[Z]}{(f_p(Z,1))} \cong F_{1,p}. $$
  Here $f_p(Z,W) \in \Z[Z,W]$ is the polynomial defined in Lemma \ref{cubesumpsin}.
\end{theorem}

\begin{proof}
    The proof of the first part is same as in Theorem \ref{2.3}. \\
    Let $P_A:=(x(P_A), y(P_A)) \in E_A[p]$ be a nontrivial element. We note that 
    $y(P_A)^2$ is a root of the polynomial $f_p(Y,A)$ which is a polynomial in $Y$ of degree $\frac{p^2-1}{6}$. 
    Using an argument similar to the proof of Theorem \ref{2.3}, we have $[K(y(P_A)^2):K]= \frac{p^2-1}{6}$,  which proves the irreducibility of $f_p(Y,A)$ over $K$. Similarly, we also conclude that $f_p(X^3+A,A)$ is irreducible over $K$.

    Since $f_p(Y,A)$ is a homogeneous polynomial with $y(P_A)^2$ being a root, we have that $\frac{y(P_A)^2}{A}$ is a root of $f_p(Y,1)$, and hence $F_{A,p}:= K(y(P_A)^2)\cong \frac{K[Z]}{(f_p(Z,1))} \cong F_{1,p}$. Since $[K_{A,p} : K]=\frac{p^2-1}{3}$, and $y(P_A)$ is a root of $f_p(Y^2,A)$, we deduce that $f_p(Y^2,A)$ is irreducible, and hence $K_{A,p} \cong \frac{K[Y]}{(f_p(Y^2,A))} := K(y(P_A))$. Similarly, since $x(P_A)$ is a root of $f_p(X^3+A,A)$ and $[K(x(P_A)):K]=\frac{p^2-1}{2}$, we have $K^\prime_{A,p} \cong \frac{K[X]}{(f_p(X^3+A,A))} := K(x(P_A))$ and $F_{A,p}=K_{A,p} \cap K^\prime_{A,p}$. \qedhere
    
\end{proof}

Next we consider the case $p \equiv 1 \pmod 6$ and let $p\OO_K=\p \overline{\p}$ where $\p$ and $\overline{\p}$ are the primes in $\OO_K$ with $\p \nmid 6A$ lying above $p$. We can write $\p= (\pi)$ where $\pi \in \Z[\zeta_3]$ such that $\pi \equiv 2 \pmod 3$. 
We define $\pi_A = -\overline{\big( \frac{4A}{\pi}\big)_6} \pi$, where $\big( \frac{\cdot}{\pi} \big)_6$ denotes the $6^{th}$ power residue symbol. Now the Gr\"ossencharacter associated to $E_A$ is given by $\psi_{E_A/K}(\p)= \pi_A$(resp. $\psi_{E_A/K}(\overline{\p})=\overline{\pi_A})$. Writing $p= a^2+ ab+ b^2$ with $a,b \in \mathbb{N}$ and $a$ odd, we see that $\pi =  \zeta_6^k( a - b \zeta_3^{i})$ for a particular choice of $0 \le k \le 5$ and $1 \le i \le 2$, hence $\pi_A =  \zeta_6^t( a - b \zeta_3^i) $ for some $0 \le t \le 5$. Define
$$ \widetilde{h_{A,p}} = \begin{cases}
\frac{\psi_{2b} \psi_a^4 - \psi_{2a} \psi_b^4}{\psi_a \psi_b} = \frac{\psi_{2b}}{\psi_{b}}\psi_{a}^3 - \frac{\psi_{2a}}{\psi_2 \psi_{a}} (\frac{\psi_{b}}{\psi_2})^3 (\psi_2^2)^2 & \text{ if } b-a \text{ is odd,}\\
\frac{\psi_{2b} \psi_a^4 - \psi_{2a} \psi_b^4}{\psi_a \psi_b \psi_2} = \frac{\psi_{2b}}{\psi_2\psi_{b}}\psi_{a}^3 - \frac{\psi_{2a}}{\psi_2 \psi_{a}} (\psi_{b})^3  & \text{ if } b-a \text{ is even, }
\end{cases}
$$
here $\psi_n$ denotes the $n$-th division polynomial for the elliptic curve $E_A$. From Lemma \ref{cubesumpsin}, it follows that $\widetilde{h_{A,p}}$ is a polynomial with integer co-efficient in $Y$ of degree $t_p=a^2+b^2 -1$ (resp. $t_p=a^2+b^2 -2$) if $b$ is even (resp. odd) and moreover, $\widetilde{h_{A,p}}(Y) = h_{A,p}(Y^2) \in \Z[Y^2]$ . In-particular, we have $\frac{2(p-1)}{3} \le t_p < p-1$, and the equality holds if and only if $|a-b| \in \{ 1, 2 \}$.

\begin{theorem}\label{3.3}
    Let $p \equiv 1 \pmod 6$ be a prime, co-prime to $A$. Let $K=\Q(\zeta_3)$ and $\p$ be a prime in $\OO_K$ lying above $p$. Let $L_{A, \p}:= K(E_A[\pi_A])$ and $F_{A, \p}$ be the unique subfield of $L_{A,\p}$ which contains $K$ and satisfies $\Gal(L_{A,\p}/F_{A,\p}) \cong \frac{\Z}{6\Z}$.  Let $\chi_A$ be the character which gives the action of $G:= \Gal(L_{A,\p}/K) \cong \frac{\Z}{(p-1)\Z}$ on $E_A[\pi_A] \cong \F_p$ and $r_\p(E_A):= \dim_{\F_p} Cl(L_{A,\p}/F_{A, \p})[p](\chi_A)$. Then we have
    $$r_\p(E_A) \le s_p(E_A) := \dim_{F_p} {\mathrm{Sel}}^p(E_A/\Q) \le 1+ r_\p(E_A),$$
    and hence $\omega(E_A)$ and $r_\p(E_A)$ uniquely determine $s_p(E_A)$.

     Furthermore, let $K_{A,\p}$ (resp. $K^\prime_{A,\p}$)  be the subfield of $L_{A,\p}$ such that $ F_{A,\p} \subset K_{A,\p} \subset L_{A,\p}$ (resp. $ F_{A,\p} \subset K^\prime_{A,\p} \subset L_{A,\p}$)  with $\Gal(L_{A,\p}/K_{A,\p}) \cong \Z/{3\Z}$ (resp. $\Gal(L_{A,\p}/K^\prime_{A,\p}) \cong \Z/{2\Z}$), then
  $$ K_{A,\p} \cong \frac{\Q[Y]}{(f_{A,p}(Y^2))}, K^\prime_{A,\p} \cong \frac{\Q[X]}{(f_{A,p}(X^3+A))} \text{ and }  F_{A,\p} = K_{A,\p} \cap K^\prime_{A,\p} \cong \frac{\Q[Y]}{(f_{1,p}(Y))} \cong F_{1,\p}, $$
 where $f_{A,p}(Y)$ is the irreducible factor (over $\Q$) of degree $\frac{p-1}{3}$ of the polynomial $h_{A,p}(Y) \in \Z[Y]$.
    

    
\end{theorem}

\begin{proof}
     The proof of the first part is same as in Theorem \ref{2.3}. \\
     Now let $P_A=(x(P_A), y(P_A))$ be a generator of $E_A[\pi_A]$. Since $\pi_AP_A=0$, it follows that $y(aP_A) =  y(bP_A)$ and hence from Lemma \ref{cubesumpsin}, it follows that $y(P_A)$ is a root of the polynomial $\widetilde{h_{A,p}}(Y)= h_{A,p}(Y^2)$.\\
     Note that $[K(x(P_A),y(P_A)): K(y(P_A))] \le 3$ and since $y(P_A)$ is root of a polynomial of degree $t_p < p-1$, it follows that $K_{A,\p} = K(y(P_A))$ and $h_{A,p}(Y^2)$ has an irreducible factor of degree $\frac{p-1}{3}$ over $K$. Thus one of the following three situations hold: $(a)$ $h_{A,p}(Y^2)$ has an unique irreducible factor of degree $\frac{p-1}{3}$ over $\Q$, $(b)$ $h_{A,p}(Y^2)$ has two distinct irreducible factors of degree $\frac{p-1}{3}$ over $\Q$,  or $(c)$ $h_{A,p}(Y^2)$ has an unique irreducible factor of degree $\frac{2(p-1)}{3}$ over $\Q$ which factors as $h(Y) \overline{h(Y)}$ over $K$. Considering $Q_A=\left(x(Q_A), y(Q_A)\right)$ a generator of $E_A[\overline{\pi_A}]$ as in Theorem \ref{2.4}, and proceeding as in the proof of Theorem \ref{2.4}, we conclude that $h_{A,p}(Y^2)$ has an unique irreducible factor $g_{A,p}(Y)$ of degree $\frac{2(p-1)}{3}$ over $\Q$ which factors as $g_{A,p}(Y)=h(Y) \overline{h(Y)}$ over $K$. \\
     Since the minimal polynomial of $y(P_A)$ has degree $\frac{2(p-1)}{3}$ over $\Q$, it follows that the minimal polynomial of $y(P_A)^2$ (over $\Q$) has degree either $\frac{2(p-1)}{3}$ or $\frac{p-1}{3}$. Since $y(P_A)^2$ is a root of $h_{A,p}(Y)$ which has degree less than $\frac{p-1}{2}$, it follows that the minimal polynomial of $y(P_A)^2$ (over $\Q$) has degree $\frac{p-1}{3}$. Consequently, we see that $h_{A,p}(Y)$ has an irreducible factor $f_{A,p}(Y)$ of degree $\frac{p-1}{3}$ (which is necessarily unique) and $g_{A,p}(Y)= f_{A,p}(Y^2)$.\\
     Since $h_{A,p}(Y)$ is a homogeneous polynomial in $Y$ and $A$, it follows that $f_{A,p}(Y)$ is a homogeneous polynomial in $Y$ and $A$. Consequently, $F_{1,\p} = K(y(P_1)^2) \cong \frac{\Q[Y]}{(f_{1,p}(Y))} = K(\frac{y(P_A)^2}{A}) = K\left(y(P_A)^2\right) :=F_{A, \p}$ and $K_{A,\p}= K(y(P_A)) \cong \frac{\Q[Y]}{(f_{A,p}(Y^2))}$. Finally, since $y(P_A)$ is a root of $f_{A,p}(Y^2)$, it follows that $x(P_A)$ is a root of $f_{A,p}(X^3+A)$ and we have $K^\prime_{A,\p}= K(x(P_A)) \cong \frac{\Q[X]}{(f_{A,p}(X^3+A))}$ , $F_{A,\p} = K_{A,\p} \cap K^\prime_{A,\p}$. \qedhere

\end{proof}

As a consequence of Theorem \ref{3.2} and \ref{3.3}, we derive the following results. The proof follows the same argument as in corollaries \ref{2.5} and \ref{2.6}.

\begin{corollary}\label{3.4}
    Suppose that $\omega(E_A) = +1$. Then $L(E_A, 1) \neq 0$ $\iff$ $\rk_\Q E_A =0$ and $\Sh(E_A/\Q)$ is finite   $\iff$ for all but finitely many primes $\p$ in $K$ with $\p \nmid 6A$, $r_\p(E_A) =0$ $\iff$ there exists a prime $\p$ in $K$ with $\p \nmid 6A$ with $r_\p(E_A) =0$ .
\end{corollary}

\begin{corollary}\label{3.5}
    Suppose that $\omega(E_A) = -1$. Then, $L(E_A, 1) = 0$ and $L^\prime(E_A,1) \neq 0$  $\iff$ $\rk_\Q E_A =1$ and $\Sh(E_A/\Q)$ is finite    $\iff$  for all but finitely many primes $\p$ in $K$ with $\p \nmid 6A$ with $r_\p(E_A) \le 1$ $\iff$ there exists a prime $\p \nmid 6A$ in $K$ lying over $p \equiv 1 \pmod 6$  with $r_\p(E_A) \le 1$.  
\end{corollary}

Observe that for every $A$ and a prime $\p \nmid 6A$ in $K$, we have the following fields diagram:
\begin{figure}[H]
\centerline{
    \xymatrix{
 & L_{A,\p}  & & & \\
 K^\prime_{A,\p} \ar@{-}^{2}[ur] & & & K_{A,\p}\ar@{-}_{3}[ull]\\
 & & F_{A, \p}\ar@{-}[ull]\ar@{-}[ur] \\
 & & K \ar@{-}[u]}}
 \caption{}
 \label{fig1}
\end{figure} 

Further note that if $A_1 = A\alpha^4$ for some $\alpha \in K$, $L_{A_1,\p} \cong L_{A,\p}$, if $A_2 = A \beta^2$ for some $\beta \in K$, then $K_{A_2,\p} \cong K_{A,\p}$ and if $A_3 = A \gamma^3$ for some $\gamma \in K$, then $K^\prime_{A_3,\p} \cong K^\prime_{A,\p}$. Consequently the following result shows that for quadratic twists (resp. cubic twists) of $E_A$, the variation of $p$-Selmer rank is controlled by a relative quadratic extension (resp. a relative cubic extension) and the root number.

\begin{proposition}\label{3.6}
      Let $\p \nmid 6A$ be a prime in $K$ lying over $p$. 
     \begin{enumerate}
         \item Let $D$ be a square free integer co-prime to $6A$ and $E_A^D: Y^2=X^3+AD^3$ denote the quadratic twist of $E_A$ by $D$. Then we have,
    $r_\p(E_A^D) \le s_p(E_A^D) \le 1 + r_\p(E_A^D)$  where $r_\p(E_A^D)= \dim_{\mathcal{O}_K/ \p} Cl(L_{AD^3,\p}/K^\prime_{AD^3,\p})[p](\chi_{AD^3})$ with $K^\prime_{AD^3,\p} \cong K^\prime_{A,\p}$. 
         \item Let $M$ be a cube free integer co-prime to $6A$ and $E_A^M: Y^2=X^3+AM^2$ denote the cubic twist of $E_A$ by $M$. Then we have,
    $r_\p(E_A^M) \le s_p(E_A^M) \le 1 + r_\p(E_A^M)$  where $r_\p(E_A^M)= \dim_{\mathcal{O}_K/ \p} Cl(L_{AM^2,\p}/K_{AM^2,\p})[p](\chi_{AM^2})$ with $K_{AM^2,\p} \cong K_{A,\p}$. 
     \end{enumerate}    
\end{proposition}

We give an application of the above result for the rational cube sum problem. We say that a natural number $n$ is a cube sum (resp. not a cube sum) if the equation $X^3+Y^3=n$ has (resp. does not have) a rational solution. It is well known that a cube free natural number $n>2$ is a cube sum iff $\rk_\Q E_{-432n^2} >0$. Using Proposition \ref{3.6}(2), we relate cube sum problem to relative cubic class group.

\begin{corollary}\label{3.7}
 Suppose $\omega(E_{-432}^{n})= 1$ and there exists a prime $\p \nmid 6n$ such that $r_{\p}(E_{-432}^n)=0$, then $n$ is not a rational cube sum.\\
 In particular, if a cube free integer $n >2 $ is a cube sum for which $\omega(E_{-432}^n)=1$, then for all primes $\p \nmid 6n$, we have $p \mid |Cl(L_{-432n^2,\p}/K_{-432,\p})|$.
 \end{corollary}

In particular since there are infinitely many primes $\ell \equiv 1 \pmod 9$ which are cube sum \cite{JMS2}, we obtain

\begin{corollary}\label{3.8}
    Given any prime $p \ge 5$, there are infinitely many cyclic cubic extensions $L_\p$ of $K_{-432,\p} \cong K_{1, \p}:=K_\p$ such that $p$ divides the relative class number $L_\p/K_\p$.
\end{corollary}

\begin{corollary}\label{3.9}
    If $\omega(E_{-432}^{n})= -1$ and there exists a prime $p \equiv 1 \pmod{6}$, which is co-prime to $n$ such that $r_\p(E_{-432}^{n}) \le 1$, then $n$ is a rational cube sum.
\end{corollary}

\section{Construction of relative class groups with large $p$-rank}\label{S4}

In this section we use known families of CM elliptic curves with large rank to construct number fields with large rank $p$-class group using Theorem \ref{2.3}, Theorem  \ref{2.4}, Theorem  \ref{3.2} and Theorem  \ref{3.3}. At first, following the construction of Kihara \cite{Kihara2004}, we show that given any odd prime $p$, there are infinitely many elliptic curves $E$ with $j(E)=1728$ for which $p$-Selmer rank of $E$ is at least $7$.

\begin{proposition}\label{4.1}
    Let $p$ be an odd prime. There exists infinitely many (non-isomorphic over $\Q$) elliptic curves $E_D: Y^2 = X^3-DX$ with the following properties:
    \begin{enumerate}
        \item $E_D$ has good reduction at $p$;
        \item $\omega(E_D) = -1$;
        \item $\rk_\Q E_D \ge 6$.
    \end{enumerate}
\end{proposition}

\begin{proof}
We recall the construction of Kihara \cite[Theorem 1]{Kihara2004} of rank $6$ (at-least) elliptic curves over $\Q(r,s,t)$. Define
$$u = \frac{rs(2t^4 -r^4 -s^4)}{t(r^4-s^4)},$$
$$a = \frac{r^4+s^4-t^4-u^4}{2(r^2s^2-t^2u^2)},$$
$$b= \frac{(s^2t^2-r^2u^2)(r^2t^2-s^2u^2)}{r^2s^2-t^2u^2},$$
then the elliptic curve $\mathcal{E}: Y^2=X^3+b(a+1)(a-1)X$ has rank at-least $6$ over $\Q(r,s,t)$. A simplification in SAGE yields that $b(a+1)(a-1) =-k(r,s,t)$, where
\begin{eqnarray*}
        k(r,s,t) &=& 2^{-8}(r^4-s^4)^{-6}r^{-4}s^{-4}t^{-12}(r^4+s^4+2r^2t^2)(r^4+s^4-2r^2t^2)(r^4+s^4+2s^2t^2)\\& & \times (r^4+s^4-2s^2t^2)
    [(r^4+s^4)^2+ 4r^2s^2t^2(r^2+s^2+t^2)] \\ & & \times ~[(r^4+s^4)^2+ 4r^2s^2t^2(r^2-s^2-t^2)]
     [(r^4+s^4)^2- 4r^2s^2t^2(r^2+s^2-t^2)]\\ & & \times ~[(r^4+s^4)^2- 4r^2s^2t^2(r^2-s^2+t^2)].
\end{eqnarray*}
 Setting $D(t):=k(2p, 1, t)$, from the specialization theorem one obtains that for all but finitely many $t \in \Q$, $\rk_\Q E_{D(t)} \ge 6$. Further from Faltings theorem \cite{Faltings1983}, it follows that for a fixed $t_0 \in \Q$,  $E_{D(t)} \cong E_{D(t_0)}$ for finitely many $t \in \Q$. Thus to complete the proof, it is enough to show that there are infinitely many $t \in \Q$ for which $\nu_p(D(t)) \equiv 0 \pmod{4}$ and $\omega(E_{D(t)}) = -1$.\\
Let us denote by $c= 16p^4+1$ and $\ell_i$ denote a prime $\ell_i \equiv 3 \pmod{4}$ and $\nu_{\ell_{i}}(c-2) \equiv 1 \pmod{2}$. Since $c-2 \equiv 3 \pmod{4}$, it follows that the number of distinct primes $\ell_i$ with the above mentioned properties is odd, say $\ell_1, \dots, \ell_{2k+1}$ and let us define $d= 3 \ell_1 \cdots \ell_{2k+1}$. A simple calculation yields that
\begin{eqnarray*}
    D(4cdpt) & \doteq & (c-2)^2(1+ 8 \alpha c p^2t^2)(1- 8 \alpha c p^2t^2)(1+ 2 \alpha c t^2)(1- 2 \alpha c t^2)\\& & \times [1+ 16 \alpha(4p^2+1)p^2t^2 + 16 \alpha^2 c^2p^2t^4] 
[1- 16 \alpha(4p^2+1)p^2t^2 + 16 \alpha^2 c^2p^2t^4] \\& & \times [1+ 16 \alpha(4p^2-1)p^2t^2 - 16 \alpha^2 c^2p^2t^4][1- 16 \alpha(4p^2-1)p^2t^2 - 16 \alpha^2 c^2p^2t^4]
\end{eqnarray*}
where $\alpha = 16p^2d^2$ and by $q_1 \doteq q_2$ we understand that $\frac{q_1}{q_2} \in (\Q^\ast)^4$.
Writing $t = \frac{m+n}{n}$ and by denoting $D(m,n) = D(4cdp \frac{m+n}{n})$, we see that 
\begin{equation*}\label{dmn}
    D(m,n) \doteq (c-2)^2 f(m,n),
\end{equation*}
where
\begin{multline*}
    f(x,y) = (y^2 + 8 \alpha c p^2 (x+y)^2)(y^2 - 8 \alpha c p^2 (x+y)^2) (y^2+ 2 \alpha c (x+y)^2)(y^2- 2 \alpha c (x+y)^2)\\
    \times [y^4+ 16 \alpha(4p^2+1)p^2y^2(x+y)^2 + 16 \alpha^2 c^2p^2(x+y)^4] [y^4- 16 \alpha(4p^2+1)p^2y^2(x+y)^2 + 16 \alpha^2 c^2p^2(x+y)^4]\\
    \times [y^4+ 16 \alpha(4p^2-1)p^2y^2(x+y)^2 - 16 \alpha^2 c^2p^2(x+y)^4][y^4- 16 \alpha(4p^2-1)p^2y^2(x+y)^2 - 16 \alpha^2 c^2p^2(x+y)^4].
\end{multline*}

Since $f(x,y)$ is a homogeneous polynomial in $2$-variables such that each of the irreducible factors has degree at most $4$, from  \cite{Greaves1992}, it follows that there are infinitely many $(m,n) \in \mathbb{N}^2$ such that $f(m,n)$ is a square-free integer. Note that the fact $f(m,n)$ is square-free necessarily implies that $(n, 2pcdm)=1$, which in-turn implies that $p \nmid (c-2)^2 f(m,n)$, thus $E_{D(m,n)}$ has good reduction at $p$.

Finally, we use the formula of root number of $E_D$, given in \cite{Liverance1995}, to show that for $(m,n) \in \mathbb{N}^2$ for which $f(m,n)$ is square-free, we have $E_{D(m,n)}$ has root number $-1$. One can check that if  $x>0$ and $y > 0$, then $f(x,y) >0$ and hence,  we have $\omega_{\infty}(E_{D(m,n)})= sgn(-(c-2)^2 f(m,n)) = -1$. Further, as $(c-2)^2f(m,n) \equiv 1 \pmod{16}$, it follows that $\omega_2(E_{D(m,n)}) = -1$. The set of primes $\ell \ge 3$ such that $\ell \equiv 3 \pmod{4}$ and $ \nu_{\ell}(D(m,n)) \equiv 2 \pmod{4}$ is precisely $\{ \ell_1, \dots, \ell_{2k+1} \}$ and for each $\ell_i$, we have $\omega_{\ell_i}(E_{D(m,n)}) = -1$. For all other primes $q$, we have $\omega_{q}(E_{D(m,n)}) = 1$, consequently, we obtain $\omega(E_{D(m,n)}) = -1$.   \qedhere

\end{proof}

\begin{theorem}\label{4.2}
 \begin{enumerate}
     \item Let $p$ be a prime such that $p \equiv 3 \pmod{4}$. There are infinitely many number fields $L_{n, p}$ of degree $4$ over $F_{n, p} \cong F_{1,p}$, such that $\rk Cl(L_{n,p})[p] \ge 12 + \rk Cl(F_{1,p})[p]$, here $F_{n,p}$ is the field defined in Theorem \ref{2.3}.
     \item Let $p$ be a prime such that $p \equiv 1 \pmod{4}$. There are infinitely many number fields $L_{n, \p}$ of degree $4$ over $F_{n, \p} \cong F_{1,\p}$, such that $\rk Cl(L_{n,\p})[p] \ge 6 + \rk Cl(F_{1,\p})[p]$, here $F_{n,\p}$ is the field defined in Theorem \ref{2.4}.
 \end{enumerate}
 \end{theorem} 

\begin{proof}
    Given an odd prime $p$, let $S_p$ be the set of fourth-power free integers $D$ such that $p \nmid D$, $\rk_\Q E_D \ge 6$ and $\omega(E_D) = -1$. The set $S_p$ is an infinite set by Proposition \ref{4.1}.\\
    For $n \in S_p$, define $L_{n,\p} := K(E_n[\p])$, where $\p$ is a prime in $K:=\Q(i)$ lying over $p$. Since $L_{n,\p}$ is a cyclic extension of degree $(N_{K/\Q}(\p)-1)$ of $K$, we have  $$\mathrm{Disc}_{L_{n,\p}/\Q}= (\mathrm{Disc}_{K/\Q})^{[L_{n,\p}: K]} N_{K/\Q}(\mathrm{Disc}_{L_{n,\p}/K}),$$ and from the conductor-discriminant formula, we have $\mathrm{Disc}_{L_{n,\p}/K} = \prod\limits_{i=1}^{N_{K/\Q}(\p)-1} \mathrm{cond}(\chi^i)$. It follows that the primes that ramify in $L_{n,\p}/\Q$ are precisely the primes dividing $2np$, since $\mathrm{cond}(\chi) = f \p$, where $f$ is the conductor of the Hecke character of $K$ associated to $E_n$ (note that conductor of $E_n$ is $|\mathrm{Disc}_{K/\Q}|N_{K/\Q}(f)$). Hence the number fields $L_{n,\p}$ are non-isomorphic for $n \in S_p$. \\
    Since for $n \in S_p$, we have $\omega(E_n) = -1$ and hence $s_p(E_n)=\dim_{\F_p} \Sel_{p}(E_n/\Q)$ is odd, and moreover as $\rk_\Q E_n \ge 6$, we see that $s_p(E_n) \ge 7$. Consequently from Theorem \ref{2.3} and Theorem \ref{2.4}, it follows that $ r_{\p}(E_n) \ge s_p(E_n) -1 \ge 6$ and hence $Cl(L_{n,\p}/F_{n,\p})[p]$ contains at least $6$-copies of $E_n[\p]$. The result follows from the fact that $E_n[\p] \cong \Z/p \Z$ for $p \equiv 1 \pmod{4}$ and $E_n[\p] \cong \Z/p \Z \times \Z/p \Z$ for $p \equiv 3 \pmod{4}$. \qedhere
    
\end{proof}

In the following proposition, we make a  quantitative estimate of the lower bound of number of number fields $L_{n,\p}$ which satisfies the condition of Theorem \ref{4.2}. For simplicity of notation, for a number field $L$, by $\mathcal{D}_L$, we denote the quantity $\mathcal{D}_L = |\mathrm{Disc}_{L/\Q}|^{\frac{1}{[L:\Q]}}$. For a prime $p \equiv 1 \pmod 4$, let 
$$\mathcal{N}_p(X) = \#\{ L \supset F_\p \cong F_{1,\p} \mid 0< \mathcal{D}_L < X, \Gal(L/F_\p) \cong \frac{\Z}{4\Z}, \rk_p \, Cl(L/F_\p)[p] \ge 6    \} $$
and, for  a prime $p \equiv 3 \pmod 4$, let 
$$\mathcal{N}_p(X) = \#\{ L \supset F_\p \cong F_{1,\p} \mid 0< \mathcal{D}_L < X, \Gal(L/F_\p) \cong \frac{\Z}{4\Z}, \rk_p \, Cl(L/F_\p)[p] \ge 12    \} .$$

\begin{proposition}\label{4.3}
    With the notation as above, for any odd prime $p$, we have
    $$\mathcal{N}_p(X) \gg  c_p X^{\frac{1}{12}} ~~\text{ where } c_p \neq 0 \text{ and }~~ X \gg 0 .$$
\end{proposition}

\begin{proof}
   From \cite[Theorem (i)]{Greaves1992}, we see that the number of square-free values $f(a,b)=n$ with $|n|=|f(a,b)| \le X$ is at least $dX^{\frac{1}{12}}$, for a non-zero constant $d$. Since ${\mathrm{Disc}}_{L_{n,\p}/K} = \prod\limits_{i=1}^{N_{K/\Q}(\p)-1} \mathrm{cond}(\chi^i)$, $\mathrm{cond}(\chi) = f \p$, and, the conductor of $E_n = 4N_{K/\Q}(f) \leq 2^6n^2$, it follows that $|\mathrm{Disc}_{L_{n,\p}/\Q}| \le (4 N_{K/\Q}(f) N_{K/\Q}(\p))^{N_{K/\Q}(\p)-1} \le (8n \sqrt{N_{K/\Q}(\p)})^{2(N_{K/\Q}(\p)-1)}$.
   Consequently, we have that $\mathcal{D}_{L_{n,\p}} \le 8|n|\sqrt{N_{K/\Q}(\p)}$ and the result follows.  \qedhere
\end{proof}

\begin{rem}
For a number field $F$, let $H_0 =F$ and inductively define $H_n$ as the $p$-Hilbert class field of $H_{n-1}$ and $H_\infty = \cup_{n \ge 0} H_n$, the maximal unramified $p$-extension of $F$. By work of \cite{GS, Boston, Hajir}, it is known that if $\rk Cl(F)[p] \ge 2 + 2 \sqrt{1+ \rk (\mathcal{O}_F^\ast)}$, then $\Gamma:= \Gal(H_\infty/F)$ is infinite, and moreover,  $\Gamma$ is not $p$-adic analytic and $\rk Cl(H_n)[p]$ is unbounded as $n$ tends to infinity. Unfortunately, it is not known if there are infinitely many number fields $F$ of (fixed) degree $d$ (which is co-prime to $p$) for which the maximal unramified $p$-extension is infinite.

By applying Theorem \ref{4.2}(2) (resp. Theorem \ref{4.2}(1)) to $p=3$ (resp. $p=5$), we see that there are infinitely many number fields $F$ of degree $16$ (resp. of degree $8$) containing $\Q(\zeta_{12})$ (resp. containing $\Q(i)$) whose maximal unramified $3$-extension (resp. maximal unramified $5$-extension)  $H_{\infty}$ is an infinite extension of $F$. Further, in this situation, $\Gamma=\Gal(H_\infty/F)$ is not $3$-adic (resp. $5$-adic) analytic group and $\lim_{n \to \infty} \rk Cl(H_n)[3] = \infty$ (resp. $\lim_{n \to \infty} \rk Cl(H_n)[5] = \infty$). 
\end{rem}

We can obtain similar results using the large rank Mordell curves constructed by Mestre.

\begin{proposition}\label{4.4}
    Let $p \ge 5$ be an odd prime. There exists infinitely many (non-isomorphic over $\Q$) elliptic curves $E_A: Y^2 = X^3+A$ with the following properties:
    \begin{enumerate}
        \item $E_A$ has good reduction at $p$;
        \item $\rk_\Q E_A \ge 7$.
    \end{enumerate}
\end{proposition}

\begin{proof}
    We recall the construction of Mestre \cite[Theorem 2]{mestre1992} of rank at-least $7$ elliptic curves over $\Q(v)$ given by
    \begin{equation*}\label{ifcs}
    \mathcal{E}: Y^2 = X^3 + (v^6-1)^2 Q(v^3)^2 P(v^3):= X^3 + A(v),
     \end{equation*}
    where $Q(v^3)$ and $P(v^3)$ are irreducible polynomials in $\Z[v]$ with $Q(v) = v^9 - v^8 - 8v^7 - 16v^6 + 194v^5 - 26v^4 + 336v^3 + 40v^2 - 11v + 3$ and $P(v)= 3^27^2v^{156} + 2^65^3v^{155} + \cdots + 3^{12}5^2$. \\
    By the specialization theorem of Silverman, we have that for all but finitely many $v_0 \in \Q$, $\rk_\Q E_{A(v_0)} \ge 7$. Further, from Faltings theorem,  one sees that given a $v_0 \in \Q$, we have $E_{A(v)} \cong E_{A(v_0)}$ for finitely many $v \in \Q$. \\
    For $p=5$, we choose $v = \frac{1}{5k}$ as  $k$ varies over the set of natural numbers, which gives us the required family of elliptic curves. For $ p >5$, we choose $v = pk$ as  $k$ varies over the set of natural numbers and obtain the required family of elliptic curves. \qedhere
    
\end{proof}

As a consequence of Proposition \ref{4.4}, one proves (in a similar fashion to that of Theorem \ref{4.2}):

\begin{theorem}\label{4.5}
 \begin{enumerate}
     \item Let $p$ be an odd prime such that $p \equiv 5 \pmod{6}$. There are infinitely many number fields $L_{n,p}$ of degree $6$ over $F_{n,p} \cong F_{1,p}$, such that $\rk Cl(L_{n,p})[p] \ge 12 + \rk Cl(F_{1,p})[p]$, where $F_{n,p}$ is the field defined in Theorem \ref{3.2}.
     \item Let $p$ be an odd prime such that $p \equiv 1 \pmod{6}$. There are infinitely many number fields $L_{n,\p}$ of degree $6$ over $F_{n,\p} \cong F_{1,\p}$, such that $\rk Cl(L_{n,\p})[p] \ge 6 + \rk Cl(F_{1,\p})[p]$, where $F_{n,\p}$ is the field defined in Theorem \ref{3.3}.
\end{enumerate}
    
\end{theorem}

\begin{rem}
    Assuming ABC conjecture, it is possible to give a quantitative lower bound for the number of number fields which satisfy the conditions of Theorem \ref{4.5}. More precisely, for a prime $p \equiv 1 \pmod 6$, define 
$\mathcal{M}_p(X) = \#\{ L \supset F_\p \cong F_{1,\p} \mid 0< \mathcal{D}_L < X, \Gal(L/F_\p) \cong \frac{\Z}{6\Z}, \rk_p \, Cl(L/F_\p)[p] \ge 6 \} $ and for  a prime $p \equiv 5 \pmod 6$, define
$\mathcal{M}_p(X) = \#\{ L \supset F_\p \cong F_{1,\p} \mid 0< \mathcal{D}_L < X, \Gal(L/F_\p) \cong \frac{\Z}{6\Z}, \rk_p \, Cl(L/F_\p)[p] \ge 12 \}$. Then using \cite[Theorem 3.5]{poonen2003} and following the proof of Proposition \ref{4.3}, one can show that for $X\gg 0$, we have $\mathcal{M}_p(X) \gg c'_pX^{\frac{1}{534}}$ for some non-zero constant $c'_p$. 
\end{rem}

We end the section by proving a conditional result on unboundedness of class group in degree $4$ and $6$ extensions.

\begin{proposition}
    \begin{enumerate}
        \item Fix an odd prime $p$ and assume that $\dim_{\F_p} \mathrm{Sel}^p(E_D/\Q)$ is unbounded as $D$ varies in the set of integers co-prime to $p$, where $E_D$ denotes the elliptic curve $y^2=x^3-Dx$. Let $F$ be a number field containing $F_{E_{1,\p}}$ of degree prime to $p$. Then $\rk \, Cl(L)[p]$ is unbounded as $L$ varies in the extensions of $F$ of degree $4$.
       \item Fix a prime $p \ge 5$ and assume that $\dim_{\F_p} \mathrm{Sel}^p(E_A/\Q)$ is unbounded as $A$ varies in the set of integers co-prime to $p$, where $E_A$ denotes the elliptic curve $y^2=x^3+A$. Let $F$ be a number field containing $F_{E_{1,\p}}$  of degree prime to $p$. Then $\rk \, Cl(L)[p]$ is unbounded as $L$ varies in the extensions of $F$ of degree $6$.
    \end{enumerate}
\end{proposition}

\begin{proof}
We'll prove the result for degree $4$ extension. The proof for degree $6$ extension is similar. By our assumption $s_p(E_D)$ is unbounded and hence by Theorem \ref{2.3}-\ref{2.4}, $r_p(E_D)$ is unbounded, and hence $\rk \, Cl(L_{E_{D,\p}})[p]$ is unbounded. This proves the result for $F=F_{E_{1, \p}}$. If $F \supset F_{E_{1, \p}}$, then for each $D$, let $L_D$ be a field containing $FL_{E_{D,\p}}$ such that $[L_D:F]=4$. Since $[FL_{E_{D,\p}} : L_{E_{D,\p}}]$ is co-prime to $p$, the norm map $N_{FL_{E_{D,\p}}/L_{E_{D,\p}}}: Cl(FL_{E_{D,\p}})_p \to CL(L_{E_{D,\p}})_p$ is surjective, unboundedness of $\rk \, Cl(L_D)[p]$ follows. \qedhere
\end{proof}


\section{Numerical Computation}\label{nc}

For a prime $p$ and an integer $D$ co-prime to $p$, using SAGE \cite{sage}, we compute $h_p(D):=\dim_{\F_p} Cl(L_{D,\p}/F_{D, \p})[p]$, here $\p$ is a prime in $\Z[i]$ lying over $p$. Since $h_p(D) \ge r_\p(E_D)$, it follows that if $h_p(D) = 0$ (resp. $h_p(D) \le 1$), then $r_{\p}(E_D) =0$ (resp. $r_{\p}(E_D) \le 1$), and hence the rank (resp. rank bound) of $E_D$ is computed using Corollary \ref{2.5} and Corollary \ref{2.6} (resp. Theorem \ref{2.3} and Theorem \ref{2.4}). If the computation of $h_p(D)$ was done assuming GRH in SAGE, we have marked the corresponding entry in the table by $^*$ and if the resulting rank (or rank bound) depends on $h_p(D)$ which was computed assuming GRH, we put a $^*$ on the rank to specify that it depends on GRH.

\begin{center}
    \underline{Table for $E_D: Y^2=X^3-DX$ with $-50 \leq D \leq 50$}
\end{center}

\begin{table}[H]
    \centering
    \renewcommand{\arraystretch}{0.9}
    \begin{tabular}{|p{0.5cm}|p{0.93cm}|p{.85cm}|p{.85cm}|p{1.0cm}|p{1.2cm}|p{0.3cm}|p{1.19cm}|p{1.15cm}|p{1.15cm}|p{1.30cm}|p{1.26cm}|}
\hline
 $D$ & $\omega(E_D)$& $h_3(D)$  & $h_5(D)$ & $h_{13}(D)$ &$\rk_{\Q}\, E_D$& & $\omega(E_{-D})$ & $h_3(-D)$ & $h_5(-D)$ & $h_{13}(-D)$ &$\rk_\Q \, E_{-D}$ \\
\hline
1  & 1 & 0 & 0 & 0 & 0 &  & 1 & 0 & 0 & 0 & 0 \\
\hline
2  & -1 & 0 & 0 & 0 & 1 &  & 1 & 0 & 0 & 0 & 0 \\
\hline
3  & 1 & NA & 0 & 0 & 0 &  & -1 & NA & 0 & 0 & 1 \\
\hline
5  & -1 & 0 & NA & 0 & 1 &  & -1 & 2 & NA & 0 & 1 \\
\hline
6  & -1 & NA & 0 & 0 & 1 &  & 1 & NA & 0 & 0 & 0 \\
\hline
7  & -1 & 0 & 0 & 0 & 1 &  & 1 & 2 & 0 & 0 & 0 \\
\hline
9  & 1 & NA & 0 & 0 & 0 &  & -1  & NA & 0 & 0 & 1 \\
\hline
10 & -1 & 2 & NA & 0 & 1 &  & 1 & 0 & NA & 0 & 0 \\
\hline
11 & 1 & 0 & 0 & 0 & 0 &  & 1 & 2 & 0 & 0 & 0 \\
\hline
13 & 1 & 2 & 0 & NA & 0 &  & -1  & 0 & 0 & NA & 1 \\
\hline
14 & -1 & 0 & 0 & 0 & 1 &  & 1 & 2 & 1 & 1 & $\leq 2$  \\
\hline
15 & -1 & NA & NA & 0  & 1 &  & -1 & NA & NA & 0 & 1 \\
\hline
17 & 1  & 2 & 1 & 1* & $\leq 2$ &  & 1 & 0 & 0 & 0 & 0 \\
\hline
18  & 1 & NA & 0 & 0 & 0 &  & -1 & NA & 1 & 0 & 1 \\
\hline
 19 & 1 & 0 & 0 & 0 & 0 &  & -1 & 0 & 0 & 0 & 1 \\
 \hline
21 & -1 & NA & 0 & 0* & 1 &  & -1  & NA & 0 & 0 & 1 \\
\hline
22 & -1 & 0 & 0 & 0 & 1 &  & 1 & 0 & 0 & 0 & 0 \\
\hline
23 & -1 & 4 & 0 & 0* & 1 &  & 1 & 0 & 0 & 0 & 0 \\
\hline
25 & -1 & 0 & NA & 0 & 1 &  & 1 & 0 & NA & 0 & 0 \\
\hline
26 & -1 & 2 & 0 & NA & 1 &  & 1 & 2 & 0 & NA & 0 \\
\hline
27 & 1 & NA & 0 & 0 & 0 &  & 1 & NA & 0 & 0 & 0 \\
\hline
29 & 1 & 0 & 0 & 0  & 0 &  & -1 & 0* & 0 & 0 & 1 \\
\hline
30 & -1 & NA & NA & 0* & 1* &  & 1 & NA & NA & 0 & 0 \\
\hline
31 & -1 & 0 & 0 & 0 & 1 &  & -1 & 2* & 0 & 0 & 1 \\
\hline
33 & 1 & NA & 0 & 0 & 0 &  & 1 & NA & 1 & 1 & $\leq 2$ \\
\hline
34 & -1 & 4* & 0 & 1* & 1 &  & 1  & 2* & 1 & 1* & $\leq 2$ \\ 
\hline
35 & 1 & 2* & NA & 0 & 0 &  & -1 & 0* & NA & 0*  & 1* \\
\hline
37 & -1 & 0* & 1 & 0 & 1 &  & -1  & 2* & 0 & 0* & 1 \\
\hline
38 & -1 & 2* & 0 & 0  & 1 &  & 1  & 2* & 0 & 0* & 0 \\
\hline
39 & -1 & NA & 0 & NA & 1 &  & 1 & NA & 1 & NA & $\le 2$ \\
\hline
41 & -1 & 0* & 0 & 0 & 1 &  & 1 & 0* & 0 & 0* & 0 \\
\hline

42 & -1 & NA & 0  & 0* & 1 &  &  1& NA & 0 & 0* & 0 \\
\hline
43 & 1 & 2* & 0 & 0* & 0 &  & 1 & 0* & 0 & 0* & 0 \\
\hline
45 & -1  & NA & NA & 0* & 1* &  & 1 & NA & NA & 0* & 0* \\
\hline
46 & -1 & 0* & 0 & 1* & 1 &  & 1 & 2* & 1 & 1* & $\leq 2$  \\
\hline
47 & -1  & 0* & 0 & 0* & 1 &  & -1 & 0* & 0 & 0* & 1 \\
\hline
49 & -1 & 0* & 0 & 1*  & 1 &  & -1  & 0* & 0 & 0* & 1 \\
\hline
50 & -1  & 0* & NA & 0* & 1*&  & 1 & 2* & NA & 0* & 0* \\
\hline
\end{tabular}
\end{table}

\newpage
\begin{small}
 Using SAGE, we compute $h_p(A):=\dim_{\F_p} Cl(L_{A,\p}/F_{A, \p})[p]$ (provided $p \nmid A$). Since $h_p(A) \ge r_\p(E_A)$, the rank (resp. rank bound) of $E_A$ is computed using Corollary \ref{3.4} and Corollary \ref{3.5} (resp. Theorem \ref{3.2} and Theorem \ref{3.3}). Outcomes dependent on GRH are marked with $^*$.   
\end{small}
\begin{center}
    \underline{Table for $E_A: Y^2=X^3+A$ with $-50 \leq A \leq 50$}
\end{center}
\begin{table}[H]
    \centering
    \renewcommand{\arraystretch}{0.7}
    \begin{tabular}{|p{0.5cm}|p{0.95cm}|p{.85cm}|p{.95cm}|p{1cm}|p{0.3cm}|p{1.15cm}|p{1.19cm}|p{1.25cm}|p{1.3cm}|}
\hline
 $A$ & $\omega(E_A)$& $h_7(A)$  & $h_{13}(A)$ &$\rk_\Q \, E_A$& & $\omega(E_{-A})$ & $h_7(-A)$ & $h_{13}(-A)$ & $\rk_\Q \, E_{-A}$ \\
\hline
1  & 1 & 0 & 0 &  0 &  & 1 & 0 & 0 &  0 \\
\hline
2  & -1 & 0 & 0* & 1 &  & -1 & 0 & 0* &  1 \\
\hline

3  & -1  & 1 & 0* & 1 &  & 1 & 0 & 0* &  0 \\
\hline
4  & 1 & 0 & 0 & 0 &  & -1 & 0 & 0* &  1 \\
\hline
5  & -1 & 0 & 0* & 1 &  & 1 & 0 & 0* & 0 \\
\hline

6  & 1 & 0 & 0* &  0 &  & 1 & 0 & 0* &  0 \\
\hline
7  & 1 & NA & 0* &0* &  & -1 & NA & 0* &  1* \\
\hline
8  & -1 & 0 & 0* & 1 &  & 1 & 0 & 0* & 0 \\
\hline
9  & -1 & 0 & 0* & 1 &  & 1 & 0 & 0* &  0 \\
\hline
10 & -1  & 0 & 0* &  1 &  & 1 & 0 & 0* &  0 \\
\hline

11 & -1 & 0 & 0* & 1 &  & 1 & 1* & 1* & $\le 2^*$  \\
\hline
12 & -1 & 1 & 0* & 1 &  & 1 & 0 & 0* &  0 \\
\hline
13 & 1 & 0 & NA &  0 &  & -1  & 0* & NA & 1* \\
\hline
14 & 1 & NA & 0* &  0* &  & 1 & NA & 0* &  0* \\
\hline
15 & 1 & 1 & 1* &  $\leq 2$ &  & -1  & 0* & 0* & 1* \\
\hline

16 & 1 & 0 & 0 & 0 &  & 1  & 0 & 0* & 0 \\
\hline
17 & 1 & 1 & 1* &  $\leq 2$ &  & 1 & 0* & 0* &  0* \\
\hline
18 & -1 & 0 & 0* & 1 &  & -1  & 0* & 0* &  1* \\
\hline
19 & -1 & 0* & 0* &  1* &  & -1 & 0*  & 0* &  1* \\
\hline
20 & 1  & 0 & 0* & 0 &  & -1  & 0* & 0* &  1* \\
\hline

21 & 1 & NA & 0* &  0* &  & -1  & NA & 0* & 1* \\
\hline
22 & -1 & 0 & 0* & 1 &  & -1  & 0* & 0* & 1* \\
\hline
23 & 1 & 1* & 0* & 0* &  & -1  & 1* & 0* & 1* \\
\hline
24 & 1 & 1 & 1* & $\le 2^*$  &  & 1 & 0* & 0* &  0*  \\
\hline
25 & 1 & 0 & 0* & 0 &  & -1  & 0* & 0* &  1* \\
\hline

26 & -1 & 1* & NA  & 1* &  & 1 & 1* & NA &  $\le 2^*$ \\
\hline
28 & -1 & NA & 0* & 1*&  & -1  & NA & 0* &  1* \\
\hline
29 & 1 &  0* & 0* & 0* &  & -1  & 0* & 0* & 1* \\
\hline
30 & -1 & 0* & 0* &1* &  & -1  & 0* & 1* &1* \\
\hline
31 & -1 & 0* & 0* & 1* &  & 1 & 0* & 0* & 0* \\
\hline
32 & 1 & 0 & 0* & 0* &  & 1 & 0 & 0* & 0* \\
\hline
33 & -1  & 0* & 1* & 1* &  & 1  & 0* & 0* & 0* \\
\hline
34 & 1 & 0* & 0* & 0* &  & 1 & 0* & 0* &  0* \\
\hline
35 & -1 & NA & 0* & 1* &  & -1  & NA & 0* &  1* \\
\hline

36 & -1 & 0 & 0* & 1 &  & 1 & 0 & 0* &  0 \\
\hline
37 & 1 & 1* & 1* & $\le 2^*$ &  & 1  & 0* & 0* &  0* \\
\hline
38 & -1 & 0* & 0* &  1* &  & -1  & 0* & 0* &  1* \\
\hline
39 & -1 & 0* & NA & 1* &  & 1 & 1* & NA &  $\le 2^*$  \\
\hline
40 & -1 & 0* & 0* &  1* &  & -1  & 0* & 0* &  1* \\
\hline

41 & -1 & 0* & 0* & 1* &  & 1  & 0* & 0* & 0* \\
\hline
42 & 1 & NA & 0* & 0* &  & 1 & NA & 0* & 0* \\
\hline
43 & 1 & 2* & 1* & $\leq 2^*$  &  & -1  & 1* & 0* & 1* \\
\hline
44 & -1 & 0* & 0* & 1* &  & -1  & 0* & 0* &  1* \\
\hline
45 & 1 & 0* & 0* & 0* &  & -1  & 0* & 0* & 1* \\
\hline
46 & -1 & 0* & 0* & 1* &  & 1 & 0* & 0* & 0* \\
\hline
47 & -1 & 0* & 0* & 1* &  & 1 & 1* & 1* & $\leq 2^*$  \\
\hline
48 & -1 & 0 & 0* & 1 &  & -1  & 0* & 0* & 1* \\
\hline
49 & 1 & NA & 0* &  0* &  & -1  & NA & 0* &  1* \\
 \hline
 50 & -1 & 0* & 1*  & 1* &  & -1  & 0* & 0* & 1* \\
 \hline
\end{tabular}
\end{table}

\newpage
In the following table we compute rank of congruent number elliptic curve using the $p$-rank of relative quadratic class group $Cl(L_{n^2,\p}/K_{n^2, \p})$, which we denote by $h_p(n)$. As before, $^*$ denotes that the outcome is correct provided GRH is true.

\begin{center}
    \underline{Table for $E_{n^2}: Y^2=X^3-n^2X$ with $n<50$}
\end{center}

\begin{table}[ht]
    \centering
    \begin{adjustwidth}{-1.25cm}{-0.8cm} 
        \begin{minipage}{0.45\textwidth}
            \flushleft 
            \renewcommand{\arraystretch}{0.8}
            \begin{tabular}{|c|c|c|c|c|c|}
                \hline
                $n$ & $\omega(E_{n^2})$& $h_3(n)$  & $h_5(n)$ & $h_{13}(n)$  &  $rk_\Q \,E_{n^2}$\\
                \hline
1  & 1  & 0 & 0 & 0 & 0 \\
\hline
2  & 1 & 0 & 0 & 0 & 0 \\
\hline
3  & 1 & NA & 0 & 0 & 0 \\
\hline
5  & -1 & 0 & NA & 0 & 1 \\
\hline
6  & -1 & NA & 0 & 0 & 1 \\
\hline
7  & -1 & 0 & 0 & 1 & 1 \\
\hline
10  & 1 & 0 & NA & 0 & 0  \\
\hline
11  & 1 & 0 & 0 & 0 &  0\\
\hline
13  & -1 & 0 & 0 & NA & 1 \\
\hline
14  & -1 & 0 & 0 & 0 & 1 \\
\hline
15  & -1 & NA & NA & 0 & 1 \\
\hline
17  & 1 & 0 & 0 & 0 & 0 \\
\hline
19  & 1 & 0 & 0 & 0 & 0 \\
\hline
21  & -1 & NA & 1 & 0 & 1 \\
\hline
22  & -1 & 0 & 0 & 1 & 1 \\
\hline
23  & -1 & 0 & 0 & 0 & 1 \\
\hline

\end{tabular}
        \end{minipage}
        \hfill 
        \begin{minipage}{0.45\textwidth}
            \flushright 
            \renewcommand{\arraystretch}{0.85}
           \begin{tabular}{|c|c|c|c|c|c|}
                \hline
                $n$ & $\omega(E_{n^2})$& $h_3(n)$  & $h_5(n)$ & $h_{13}(n)$  &  $rk_\Q \,E_{n^2}$\\
                \hline
                
               26  & 1 & 2 & 0 & NA & 0 \\
\hline
                29  & -1 & 2 & 0 & 0 & 1 \\
                \hline
                30  & -1 & NA & NA & 0 & 1 \\
                \hline
                31  & -1 & 2 & 0 & 0 & 1 \\
                \hline
                
                33  & 1 & NA & 0 & 0 & 0 \\
                \hline
                34  & 1 & 2 & 1 & 1 & $\leq 2$ \\
                \hline
                35  & 1 & 0 & NA & 0 & 0 \\
                \hline
                
                37  & -1 & 0 & 0 & 0 & 1 \\
                \hline
                38  & -1 & 2 & 1 & 0 & 1 \\
                \hline
                39  & -1 & NA & 0 & NA & 1 \\
                \hline
                
                41  & 1 & 2 & 1 & 1 & $\leq 2$ \\
                \hline
                42  & 1 & NA & 0 & 0 & 0 \\
                \hline
                43  & 1 & 2 & 0 & 0* & 0 \\
                \hline
                
                46  & -1 & 2 & 0 & 0* & 1 \\
                \hline
                47  & -1 & 2 & 0 & 0* & 1 \\
                \hline
                
            \end{tabular}
        \end{minipage}
        \end{adjustwidth}
       \end{table}

In the following table we compute rank of cube sum elliptic curve using the $p$-rank of relative cubic class group $Cl(L_{-432n^2,\p}/K_{-432n^2, \p})$, which we denote by $h_p(n)$.  As before, $^*$ denotes that the outcome is correct provided GRH is true.

\begin{center}
    \underline{Table for $E_{-432n^2}: Y^2=X^3-432n^2$ with $n \le 50$}
\end{center}

\begin{table}[ht]
    \centering
    \begin{adjustwidth}{-1.25cm}{-0.8cm} 
        \begin{minipage}{0.45\textwidth}
            \flushleft 
             \renewcommand{\arraystretch}{0.8}
            \begin{tabular}{|c|c|c|c|c|}
                \hline
                $n$ & $\omega(E_{-432n^2})$& $h_7(n)$  & $h_{13}(n)$  &  $rk_\Q \, E_{-432n^2}$\\
                \hline
1  & 1 & 0 & 0  & 0 \\
\hline
2  & 1 & 0 & 0  & 0 \\
\hline
3  & 1 & 0 & 0*  & 0 \\
\hline
4  & 1 & 0 & 0  & 0 \\
\hline
5  & 1 & 0 & 0*  & 0 \\
\hline
6  & -1 & 0 & 0* & 1 \\
\hline
7  & -1 & NA & 0*  & 1* \\
\hline
9  & -1 & 0 & 0*  & 1 \\
\hline
10  & 1 & 0 & 0*  & 0 \\
\hline
11  & 1 & 0 & 0*  & 0 \\
\hline
12  & -1 & 0 & 0*  & 1 \\
\hline
13  & -1 & 0 & NA  & 1 \\
\hline
14  & 1 & NA & 0*  & 0* \\
\hline
15  & -1 & 1 & 0*  & 1 \\
\hline
17  & -1 & 0 & 0* &  1 \\
\hline
18  & 1 & 0 & 0* & 0 \\
\hline
19  & 1 & 1 & 1* & $\leq 2$ \\
\hline
20  & -1 & 0 & 0* & 1 \\
\hline
21  & 1 & NA & 0* & 0* \\
\hline
22  & -1 & 0 & 0*  & 1 \\
\hline
23  & 1 & 0 & 0* & 0 \\
\hline
25  & 1 & 0 & 0*  & 0 \\
\hline
\end{tabular}
        \end{minipage}
        \hfill 
        \begin{minipage}{0.45\textwidth}
            \flushright 
             \renewcommand{\arraystretch}{0.85}
           \begin{tabular}{|c|c|c|c|c|}
                \hline
                $n$ & $\omega(E_{-432n^2})$& $h_7(n)$  & $h_{13}(n)$  &  $rk_\Q \, E_{-432n^2}$\\
                \hline
                26  & -1 & 0 & NA & 1 \\
                \hline
                28  & -1 & NA & 0* &  1* \\
                \hline
                29  & 1 & 0 & 0* &  0 \\
                \hline
                30  & 1 & 1 & 1* & $\leq 2$ \\
                \hline
                31  & -1 & 0 & 0*  & 1 \\
                \hline
                33  & -1 & 0 & 0*  & 1 \\
                \hline
                34  & -1 & 0 & 0*  & 1 \\
                \hline
                35  & -1 & NA & 0*  & 1* \\
                \hline
                36  & 1 & 0 & 0* & 0  \\
                \hline
                37  & 1 & 1 & 1*  & $\leq 2$ \\
                \hline
                38  & 1 & 0 & 0*  & 0 \\
                \hline
                39  & 1 & 0 & NA &  0 \\
                \hline
                41  & 1 & 0 & 0* &  0 \\
                \hline
                42  & -1 & NA & 0* &  1* \\
                \hline
                43  & -1 & 0 & 0* &  1 \\
                \hline
                44  & 1 & 0 & 0* &  0 \\
                \hline
                45  & 1 & 0 & 0* &  0 \\
                \hline
                46  & 1 & 0 & 0* &  0 \\
                \hline
                47  & 1 & 0 & 0* &  0 \\
                \hline
                49  & -1 & NA & 0* &  1* \\
                \hline
                50  & -1 & 0 & 0* &  1 \\
                \hline
            \end{tabular}
        \end{minipage}
        \end{adjustwidth}
       \end{table}

\newpage
\bibliographystyle{amsalpha}
\bibliography{references}

\end{document}